\documentclass[reqno,twoside]{amsart}
\usepackage{amsmath}
\usepackage{amsfonts}
\usepackage{amssymb}
\usepackage{enumerate}

\usepackage{a4wide,amsmath,amssymb,latexsym,amsthm}

\usepackage{mathrsfs,mathtools,epic,bm}
\usepackage{hyperref}
\hypersetup{colorlinks=true,linkcolor=blue,filecolor=mangeta,urlcolor= cyan}
\newtheorem{theorem}{Theorem}[section]
\newtheorem{proposition}[theorem]{Proposition}
\newtheorem{lemma}[theorem]{Lemma}
\newtheorem{corollary}[theorem]{Corollary}
\theoremstyle{definition}
\newtheorem{definition}[theorem]{Definition}

\newtheorem{remark}[theorem]{Remark}

\newtheorem{assumption}[theorem]{Assumption}

\newcommand\norm[1]{\left\lVert #1\right\rVert}
\newcommand\abs[1]{\left\lvert#1\right\rvert}

\numberwithin{equation}{section}

\def \dis {\displaystyle}

\def \R {\mathbb{R}}

\def \Ha {{H^{\alpha}_{a^+}(a, b)}}

\def \Hb {{H^{\alpha}_{b^-}(a, b)}}

\def \V {\cal V}

\def \Dc {{\mathcal D}}
\def \Dr {{\mathbb D}}
\def \D {{\mathbb D}}
\def \V {{\mathbb V}}

\def \Da { \Dr_{a^+}^\alpha}

\keywords{Riemann-Liouville fractional derivative, Caputo fractional derivative, low regret control, no-regret control, Sturm-Liouville equations, optimality system}

\begin{document}
\title[No-regret and low-regret controls in a star graph]{No-regret and low-regret controls  of space-time fractional parabolic Sturm-Liouville equations in a star graph }

\author{Gis\`{e}le Mophou}
\address{Gis\`{e}le Mophou, Laboratoire L.A.M.I.A., D\'{e}partement de Math\'{e}matiques et Informatique, Universit\'{e} des Antilles, Campus Fouillole, 97159 Pointe-\`{a}-Pitre,(FWI), Guadeloupe -
Laboratoire  MAINEGE, Universit\'e Ouaga 3S, 06 by 10347 Ouagadougou 06, Burkina Faso.}
\email[Mophou]{gisele.mophou@univ-antilles.fr}

\author{Maryse Moutamal}
\address{Maryse Moutamal, University of Buea, Department of Mathematics, Buea, Cameroon.}
\email{maryse.moutamal@aims-cameroon.org}

\author{Mahamadi Warma}
\address{Mahamadi Warma, Department of Mathematical Sciences and the Center for Mathematics and Artificial Intelligence (CMAI), George Mason University,  Fairfax, VA 22030, USA.}
\email{mwarma@gmu.edu}


\thanks{The second author is supported by the Deutscher Akademischer Austausch Dienst/German Academic Exchange Service (DAAD).  The work of the third author is partially supported by the US Army Research Office (ARO) under Award NO: W911NF-20-1-0115}

\begin{abstract}	
We are concerned with a space-time fractional parabolic initial-boundary value problem of Sturm–Liouville type in a general star graph with mixed Dirichlet and Neumann boundary controls. We first give several existence, uniqueness and regularity results of weak and very-weak solutions.  Using the notion of no-regret control introduced by Lions, we prove the existence, uniqueness, and characterize the low regret control of a quadratic boundary optimal control problem, then we prove that this low regret control converges to the no-regret control and we provide the associated optimality systems and conditions that characterize that no-regret control.
\end{abstract}

\maketitle

\section{Introduction and problem setting}

\begin{figure}[h]
	\centering
	\includegraphics[width=0.8\linewidth]{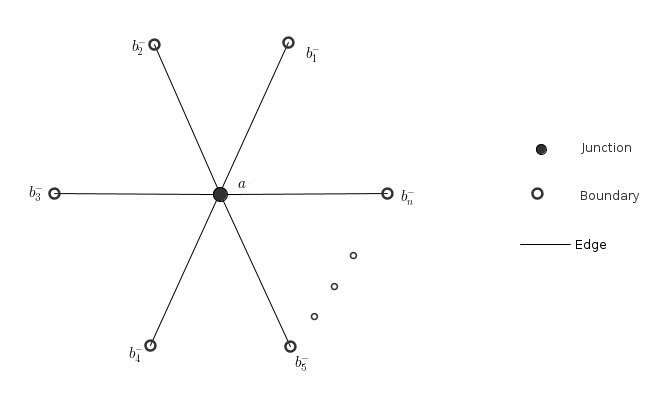}
	\caption{A sketch of a star graph with n edges}
	\label{fig:stargraph}
\end{figure}

The main concern of the present paper is to study the following optimal control problem:
\begin{equation}\label{op1}
	\dis \min_{v\in \left(L^2(0,T)\right)^{N-1}}J(v,y^0):=	\sum_{i=1}^{N}\int_{Q_i} \abs{y^i(v,y^0)-y_d^i}^2\,\mathrm{d}x\mathrm{d}t+\zeta\sum_{i=2}^{N}\int_0^T|v^i|^{2} \,\mathrm{d}t,\;\;\forall\; y^0\in \mathbb L^2,
\end{equation}
subject to the constraints that  $y(v,y^0)$ solves the space-time fractional system on a general star graph (see Figure \ref{fig:stargraph}):
\begin{equation}\label{pa7}
	\left\{
	\begin{array}{lllllllllllllllll}
		\displaystyle	\Dc_t^{\gamma}y^i+\mathcal{D}_{b_i^-}^\alpha(\beta^i\mathbb{D}_{a^+}^\alpha y^i)+q^iy^i&=&f^i&\text{in}& {(0,T)\times (a,b_i)}=:Q_i,\,i=1,\dots, N,\\
		\displaystyle	I_{a^+}^{1-\alpha}y^i(\cdot,a^+)-I_{a^+}^{1-\alpha}y^j(\cdot,a^+)&=&0&\text{in}&(0,T),~i\neq j=1,\dots, N,\\
		\displaystyle	\sum_{i=1}^N\beta^i(a)\mathbb{D}_{a^+}^\alpha y^i(\cdot,a^+)&=&0&\text{in}& (0,T), \\
		\displaystyle	I_{a^+}^{1-\alpha} y^1(\cdot,b_1^-)&=&0&\text{in}&(0,T), \\
		\displaystyle	I_{a^+}^{1-\alpha} y^i(\cdot,b_i^-)&=&v^i&\text{in}& (0,T), \; i=2,\dots, m\\
		\displaystyle	\beta^i(b_i)\mathbb{D}_{a^+}^\alpha y^i(\cdot,b_i^-)&=&v^i&\text{in}&(0,T), \; i=m+1,\dots,N,\\
		\displaystyle	y^i(0,\cdot) &=&y^{0,i}&\text{in}& (a,b_i),~~i=1,\dots,N.
	\end{array}
	\right.
\end{equation}
Here,  $N\in\mathbb N$, $0<\gamma,\alpha<1$, $\zeta >0$, the functions $f=(f^i)_i\in L^2((0,T);\mathbb{L}^2)$,  $y_d=(y_d^i)_i\in L^2((0,T); \mathbb L^2)$ are given  whereas the initial condition $y^0=(y^{0,i})_i\in \mathbb{L}^2$  is unknown. The control functions   $v^i\in L^2(0,T),\, i=2,\cdots N$ and
\begin{equation*} 
	\mathbb{L}^2:=\dis \prod_{i= 1}^N L^2(a,b_i).
\end{equation*}
In \eqref{pa7} $	\Dc_t^{\gamma}$ and $\mathbb D^\alpha$ denote the Caputo and the Riemann-Liouville fractional derivatives of order $\gamma$ and $\alpha$, respectively, and $I^\alpha$ is the Riemann-Liouville fractional  integral of order $\alpha$.

Since {$ \mathbb L^2$ } is an infinite dimensional space, the minimization problem \eqref{op1} does not always make sense.
One idea is to consider instead the following  {optimization} problem:
\begin{align}\label{op2}
	{\min_{v\in \left(L^2(0,T)\right)^{N-1}}}\left[\sup_{y^0\in  \mathbb L^2}\,J(v,y^0)\right].
\end{align}
But proceeding in that way,  it may happen that $\dis \sup_{y^0\in  \mathbb L^2}\,J(v,y^0)=+\infty$. This difficulty led J.L Lions  \cite{lions1992} to think about looking for a control $v$ such that
$$J(v,y^0)\leq J(0,y^0),\forall y^0\in \mathbb L^2.$$
Fractional derivatives provide an excellent tool for the description of memory and hereditary effects of various materials and processes \cite{Zettl,kil}. Time-fractional diffusion equations can describe anomalous diffusion on fractals see \cite{marichev,Oldham,kil} and references therein. Existence results of initial and boundary value problems for such equations have been studied widely. We refer, for instance, to \cite{Ga-Wa-2021,Oldham} and references therein.\par
The implementation of differential equations to the network domains, in particular in biology, engineering, cosmology, leads Lumer \cite{Lumer} to the notion of evolution problems on ramified spaces. Ever since, work on ordinary and partial differential equations on metric graphs has greatly evolved (see e.g \cite{vonbelow,gunter2022} and the references therein). The survey paper by D\`ager and Zuazua \cite{dager} is an excellent reference where several 1-D optimal control problems on graphs have been studied.  As an application, in \cite{mophou2018}, the authors study the control of a  gas network which  is represented by a nonlinear model derived from a semi-linear approximation of the fully nonlinear isothermal Euler gas equation. \par
Fractional optimal control problems have attracted several authors in the last two decades. Agrawal \cite{agrawal2004general} detains the first record as he proposed a general formulation and a numerical method to solve such problems.
The concepts of no-regret control and low regret control were developed by Lions \cite{lions1992,Lions2000} in order to control systems modeled by partial differential equations with integer time derivatives and uncertainties (unknown perturbations). Since then, many authors have used these concepts to control system with partial information or with missing data. For instance in \cite{Mophou2017}, Mophou investigated an optimal control of a time-fractional diffusion equation with a missing boundary condition.  Using the notion of no-regret control and low  regret control, she proved that the low regret control problem associated with the boundary fractional diffusion equation has a unique solution which converges to the no-regret control. In \cite{kenne2020optimal}, Kenne and al.  studied a model of population dynamics with age dependence and spatial structure but unknown birth rate. Using the notion of low-regret control, they proved that the state of the system can reach a desired state by acting on the system via a localized distributed control; then, using an appropriate Hilbert space, they proved that the family of low-regret controls tends to a so-called no-regret control. Nakoulima et al. \cite{Nakoulina2002,Nakoulima2003} applied these notions to control distributed nonlinear systems with incomplete data.  In all these papers, the authors proved the uniqueness of the low-regret control and its convergence to the no-regret control. Then,  they provided the optimality system that characterizes the no-regret control.\par
In this paper, we study the optimal control of a space-time fractional diffusion equation involving the Caputo time-fractional derivative and the space-fractional Sturm-Liouville operator, where the source is
missing or incomplete. We first prove the existence and uniqueness of the solution to the diffusion equation, then using the no-regret control and the low regret control notions, we first prove the uniqueness of the low regret control and we show that this control converges to the no-regret control that we characterize by a singular optimality system.\par
The rest of the paper is organized as follows. In Section \ref{prelim}, we give some preliminary results that will be used
for the proofs of the existence of solutions to fractional diffusion equations involving Sturm-Liouville derivatives and for the study of the optimal control problem. In Section \ref{oneedge}, we consider the one edge case Sturm-Liouville equation where we show existence, regularity and an explicit representation of solutions. In Section \ref{graph}, we study the existence, uniqueness, regularity and explicit representation of solutions of the Sturm-Liouville problem in a general star graph.  In section \ref{no-regret}, firstly,  we show the existence and uniqueness of the low-regret control and characterize the associated optimality systems and conditions.   Secondly, we show the existence  and uniqueness of the no-regret control and characterize the associated optimality systems and conditions.  We finish the paper with some remarks in Section \ref{conclusion}.

\section{Preliminaries}\label{prelim}

In this section, we introduce some notations,  give the function spaces needed to study our problems,  recall some known results and prove some intermediate results that are needed in the proofs of our main results. We start with fractional integrals and derivatives.
Let $[a,b]\subset\R$, $a\ge 0$ and $f : [a,b] \to \R$ be a given function.

\begin{definition}\label{defIRL}
The left,  and right Riemann–Liouville fractional integrals of order $\alpha\in(0,1)$ of $f$,  are defined, respectively,  by:
$$\begin{array}{ lll}&(I^{\alpha}_{a^+}f)(x):=\dis \frac{1}{\Gamma(\alpha)}\int_{a}^{x}(x-t)^{\alpha-1}f(t)dt,\quad (x>a)\\&(I^{\alpha}_{b^-}f)(x):=\dis \frac{1}{\Gamma(\alpha)}\int_{x}^{b}(t-x)^{\alpha-1}f(t)\mathrm{d}t,\quad (x<b),
\end{array}$$ provided that the integrals exist, and $\Gamma$ denotes the usual Euler-Gamma function.
\end{definition}

\begin{definition}\cite{kil}\label{defDRL}
The left, and right Riemann–Liouville fractional derivatives of order $\alpha\in(0,1)$ of $f$,  are defined, respectively, by
\begin{equation}\label{4}
	\begin{array}{lll}&(\Dr^{\alpha}_{a^+}f)(x):=D(I^{1-\alpha}_{a^+}f)(x)=\dis \frac{1}{\Gamma(1-\alpha)}\frac{\mathrm{d}}{\mathrm{d}x}
\int_{a}^{x}(x-t)^{-\alpha}f(t)\mathrm{d}t,\quad (x>a)
\\&
(\D^{\alpha}_{b^-}f)(x):=-D(I^{1-\alpha}_{b^-}f)(x)=\dis \frac{-1}{\Gamma(1-\alpha)}\frac{\mathrm{d}}{\mathrm{d}x}
\int_{x}^{b}(t-x)^{-\alpha}f(t)\mathrm{d}t,\quad ~~(x<b),
\end{array}
\end{equation}
	 provided that the integrals exist.
\end{definition}	
	
\begin{definition}\cite{kil}
The left, and right-sided Caputo fractional derivatives of order $\alpha\in(0,1)$ of $f$,  are defined respectively, by:
	\begin{equation}\label{4b1}\begin{array}{lll}&
(\Dc^{\alpha}_{a^+}f)(x):=(I^{1-\alpha}_{a^+}Df)(x)=\dis \frac{1}{\Gamma(1-\alpha)}\int_{a}^{x}(x-t)^{-\alpha}f'(t)\mathrm{d}t,\quad (x>a)
\\&
(\Dc^{\alpha}_{b^-}f)(x):= -(I^{1-\alpha}_{b^-}Df)(x)=\frac{-1}{\Gamma(1-\alpha)}{\dis \int_{x}^{b}}(t-x)^{-\alpha}f'(t)\mathrm{d}t\quad (x<b)\end{array}
\end{equation}
	provided that the integrals exist.
\end{definition}


\begin{definition}\label{Mittag}\cite{pod, kil}
Let $\alpha, \beta\in \mathbb{C}$ be such that $\mbox{Re}(\alpha)> 0 $ and $\mbox{Re}(\beta)>0$. The classical Mittag-Leffer function of two parameters denoted by $E_{\alpha,\beta}$ is given by
\begin{equation}\label{eq23}
	E_{\alpha,\beta}(z)=\sum_{n=0}^{\infty}\frac{z^n}{\Gamma(\alpha n+\beta)}~~~z\in\mathbb{C}.
\end{equation}
\end{definition}

\begin{lemma}\cite{pod}\label{Lmit} 
Let  $0 <\alpha < 2$, $\beta>0$ and let $\mu$ to be such that
\begin{equation*}
	\frac{\pi\alpha}{2}<\mu< min(\pi,\pi\alpha).
\end{equation*}
Then, there is a  constant $C =C(\alpha,\beta, \mu)>0$ such that
\begin{equation}\label{Mitl3}
\abs{	E_{\alpha,\beta}(z)}\leq \frac{C}{1+\abs{z}}\,\,\,\,\,\,~~\,\,\mu\leq arg(z)\leq \pi.
\end{equation}
\end{lemma}
We have the following result.  We refer to  \cite{kil,marichev} for the proof.

\begin{lemma}\label{l1}
 Let $0<\alpha<1$,  $1<p<{1}/{\alpha}$, $q={p}/{(1-\alpha p)}$,  and $\rho\in L^p (a, b)$. Then, there is a constant $C=C(\alpha,p,q,a,b)>0$ such that  for every $f\in L^p(a,b)$, 
 \begin{align*}	
\norm{I^{\alpha}_{a^+}\rho}_{L^q(a,b)}\leq&C\norm{\rho}_{L^p(a,b)},\\
\norm{I^{\alpha}_{b^-}\rho}_{L^q(a,b)}\leq&C\norm{\rho}_{L^p(a,b)}.
\end{align*}
\end{lemma}

\begin{remark}\label{rem27}
Using the Young convolution inequality we obtain that for every $0<\alpha<  1$ there is a constant $C>0$ such that for every $\rho\in L^2(a,b)$,
\begin{equation}\label{marichev13bis}
\|I^\alpha_{a^+}\rho\|_{L^2(a,b)}\leq C\|\rho\|_{L^2(a,b)},
\end{equation}
\begin{equation}\label{marichev13ter}
\|I^\alpha_{b^-}\rho\|_{L^2(a,b)}\leq C\|\rho\|_{L^2(a,b)}.
\end{equation}
\end{remark}

Next, let $c_0,d_0\in \R$. Let $f:[a,b]\to \R$ have the representation
\begin{align}\label{a1}	f(x)=\frac{c_0}{\Gamma(\alpha)}(x-a)^{\alpha-1}+I^{\alpha}_{a^+}\varphi(x)\;\mbox{ for a.e. } x\in [a,b],
	\end{align}
	 and let also $g:[a,b]\to \R$ have the representation
\begin{align}\label{a2}
	g(x)=\frac{d_0}{\Gamma(\alpha)}(b-x)^{\alpha-1}+I^{\alpha}_{b^-}\psi(x)~~\mbox{ for a.e.  }  x\in [a,b],
\end{align}
where $\varphi$ and $\psi$ belong to $L^2(a,b)$. We shall denote by $AC^{\alpha,2}_{a^+}$ and $AC^{\alpha,2}_{b^-}$ the spaces of all functions $f$ and $g$ having the representations \eqref{a1} and \eqref{a2}, respectively, with $\varphi , \psi \in L^2 (a, b)$.

\begin{remark}
We observe the following:
\begin{subequations}	
\begin{alignat}{11}
	\D^{\alpha}_{a^+}f\in L^2 (a, b)&\iff &f \in AC^{\alpha,2}_{a^+} \label{a3},\\ 	\D^{\alpha}_{b^-}f\in L^2 (a, b)&\iff& f \in AC^{\alpha,2}_{b^-} \label{a4}.
	\end{alignat}
\end{subequations}
\end{remark}

For more details  on these spaces and the proof of \eqref{a3}-\eqref{a4},  we refer to \cite{darius}.

We set
\begin{subequations}
\begin{alignat}{11}
 \Ha&=&AC^{\alpha,2}_{a^+} \cap L^2 (a, b)\label{defHa}\\
		\Hb&=&AC^{\alpha,2}_{b^-} \cap L^2 (a, b)\label{H2}.
	\end{alignat}
\end{subequations}
It follows from the definitions of $AC^{\alpha,2}_{a^+}$ and $AC^{\alpha,2}_{b^-}$ that,
\begin{subequations}
\begin{alignat}{11}
\rho \in \Ha \Longleftrightarrow \rho \in L^2(a, b) \hbox{ and } \D^{\alpha}_{a^+}\rho \in L^2(a, b),\label{prelim1}\\
	\rho \in \Hb \Longleftrightarrow \rho \in L^2(a, b) \text{ and } \D^{\alpha}_{b^-}\rho \in L^2(a, b).\label{prelim2}
\end{alignat}
\end{subequations}
We endow $\Ha$ with  the inner product
\begin{equation}\label{prodHa}
(\varphi,\psi)_{\Ha}=\int_a^b\varphi\psi \mathrm{d}x+\int_a^b\D^{\alpha}_{a^+}\varphi\D^{\alpha}_{a^+}\psi\mathrm{d}x.
\end{equation}
Then,  $\Ha$ endowed  with the norm given by
	\begin{equation}\label{normHa} \norm{\varphi}^2_{\Ha}=\norm{\varphi}^2_{L^2(a,b)}+\norm{\D^{\alpha}_{a^+}\varphi}^2_{L^2(a,b)}
	\end{equation}
is a Hilbert space  (see e.g. \cite{darius}).  Moreover,
the norm  on
	$\Ha$ given by \eqref{normHa} is equivalent to
	\begin{equation}\label{n2}		\abs{\norm{\varphi}}^2_{\Ha}=\abs{I^{1-\alpha}
\varphi(a^+)}^2+\norm{\D^{\alpha}_{a^+}\varphi}^2_{L^2(a,b)}~.
		\end{equation}

 \begin{lemma}\cite[Corollary 32]{darius}
  Let $1/2<\alpha\leq 1$ . Then, the following continuous embedding is compact
	\begin{equation}\label{a6}
	{\Ha} \hookrightarrow L^2(a,b).
	\end{equation}
\end{lemma}

We also introduce the space $\mathcal{V}$ given by
\begin{equation}\label{defVr}
\mathcal{V} = \Big\{y \in \Ha :\;  \Dc^{\alpha}_{b^-}(\beta \Dr^{\alpha}_{a^+}y) \in H_{b^-}^{1-\alpha}(a, b)\Big\},
\end{equation}
where $\beta\in C[a,b]$ and  $H_{b^-}^{1-\alpha}(a, b)$ is defined in \eqref{H2}.
We have the following trace result  that will be useful for some calculations in the upcoming sections.

\begin{lemma}\label{trace}\cite{gunter2022}
Let $0<\alpha\leq 1$.
Then, the following assertions hold.
\begin{enumerate}
\item Let $T>0$ and $\rho \in L^2((0,T);\Ha)$.  Then, for every $x_0\in [a,b],$ the function  $I_{a^+}^{1-\alpha}\rho(\cdot,x_0)$ exists   and belongs to $L^2(0,T)$. Moreover, there is a constant $C=C(a,b,\alpha)>0$ such that
  \begin{align}
    \|I_{a^+}^{1-\alpha}(\rho)(\cdot,x_0)\|^2_{L^2(0,T)}\leq C \|\rho\|^2_{L^2((0,T);\Ha)}. \label{traceI}
    \end{align}
  \item Let  $\mathcal V$ be the space defined in \eqref{defVr}.
Then, $\mathcal{V}$ endowed with the norm
\begin{align}\label{NM}
\|y\|_{\mathcal V}:=&\left(\|y\|_{\Ha}^2+\|  \Dc^{\alpha}_{b^-}(\beta \Dr^{\alpha}_{a^+}y)\|_{H_{b^-}^{1-\alpha}(a, b)}^2\right)^{1/2}\notag\\
=&\left(\|y\|_{\Ha}^2+ \|  I^{1-\alpha}_{b^-}(\beta \Dr^{\alpha}_{a^+}y)^\prime)\|_{L^2(a,b)}^2+\| (\beta\Dr_{a^+}^\alpha\rho)^\prime\|_{L^2(a, b)}^2\right)^{1/2}
\end{align}
and the associated scalar product
\begin{align}\label{SP}
(\phi,\psi)_{\mathcal V}:=&\int_a^b\phi\psi\;dx +\int_a^b \Dr_{a^+}^\alpha\phi \Dr_{a^+}^\alpha\psi\;dx +\int_a^b I^{1-\alpha}_{b^-}(\beta \Dr^{\alpha}_{a^+}\phi)^\prime I^{1-\alpha}_{b^-}(\beta \Dr^{\alpha}_{a^+}\psi)^\prime\;dx\notag\\
&+\int_a^b (\beta\Dr_{a^+}^\alpha \phi)^\prime (\beta\Dr_{a^+}^\alpha\psi)^\prime\;dx,
\end{align}
 is a Hilbert space.
 \item  Let $\rho \in \mathcal{V}$.  Then,  $\Dr_{a^+}^\alpha\rho\in C[a,b]$  and there is a constant $C=C(a,b,\alpha)>0$ such that
 \begin{equation}\label{TraceD}
 \|\Dr_{a^+}^\alpha\rho\|_{C[a,b]}\leq C\|\rho\|_{\mathcal V}.
\end{equation}
    \end{enumerate}
\end{lemma}

\begin{remark}\label{remtrace}
Let $0<\alpha\leq 1$.
Notice that if $\rho\in L^2((0,T);\Ha)$,  then  from Lemma \ref{trace} we have that the traces
$(I^{1-\alpha}_{a^+})\rho)(\cdot,a)$ and $(I^{1-\alpha}_{a^+}\rho)(\cdot,b)$
 exist for a.e. $t\in (0,T)$, and belong to $L^2(0,T)$. Also if $\rho\in \mathcal V$, then $(\beta\Dr_{a^+}^\alpha\rho)(\cdot,a)$ and $(\beta\Dr_{a^+}^\alpha\rho)(\cdot,b)$ exist and are finite.
\end{remark}

Next, we introduce some  integration by parts formulas. We refer to \cite{agrawal2004general} for the proof.

\begin{lemma}\label{lem0}
Let  $y,\phi\in
\mathcal{V}.$ Then,
 \begin{equation}\label{integ1}
 	\begin{array}{lll}
 		&\int_{a}^{b}y\D^{\alpha}_{b^-}\phi\;\mathrm{d}x+ \dis\int_a^b \phi\Dc^{\alpha}_{a^+}y\;\mathrm{d}x=\left[y(s)I^{1-\alpha}_{b^-}\phi(x)\right]_{x=a}^{x=b}
 ,
\\&
\dis \int_a^b y\Dr_{a^+}^\alpha \phi \;dx
 =\dis
 \int_a^b \phi (\Dc_{b^-}^\alpha y)\;dx+ \left[(I_{a^+}^{1-\alpha}\phi)(x)y(x)\right]_{x=a}^{x=b}.
 	\end{array}
 \end{equation}
\end{lemma}

More generally, we have the following.

\begin{lemma}\label{lem01}
Let $\beta\in {C}[a,b]$ and $y,\phi\in
\mathcal{V}.$ Then
 \begin{align}\label{integ1&2}
\dis \int_a^b \phi(s)\Dc_{b^-}^\alpha (\beta\Dr_{a^+}^\alpha y )(s) \;ds
 =&\dis
 \int_a^b (\beta\Dr_{a^+}^\alpha y )(s) \Dr_{a^+}^\alpha (\phi))(s)\;  ds-
 \left[(\beta\Dr_{a^+}^\alpha y )(s)I_{a^+}^{1-\alpha}(\phi)(s)\right]_{s=a}^{s=b}
  \notag\\
 =&\dis
 \int_a^b y(s) \Dc_{b^-}^\alpha (\beta\Dr_{a^+}^\alpha (\phi))(s)\;  ds\notag\\
& -\left[(\beta\Dr_{a^+}^\alpha y )(s)I_{a^+}^{1-\alpha}(\phi)(s)\right]_{s=a}^{s=b}
 + \dis \left[I_{a^+}^{1-\alpha}(y)(s) \beta(s)\Dr_{a^+}^\alpha (\phi)(s)\right]_{s=a}^{s=b}.
 \end{align}
 \end{lemma}

\begin{remark}\label{reminteg1}
Note that as for \eqref{integ1}, we can write
\begin{equation}\label{integ1t}
 	\begin{array}{lll}
 		&\dis\int_0^T \phi\Dc^{\gamma}_{t}y\;\mathrm{d}t=\int_{0}^{T}y\D^{\gamma}_{T^-}\phi\;\mathrm{d}t+\left[y(t)I^{1-\gamma}_{T^-}\phi(t)\right]_{t=0}^{t=T},
\\&
\dis \int_0^T y\Dr_{t}^\gamma \phi \;dt
 =\dis
 \int_0^T \phi (\Dc_{T^-}^\gamma y)\;dt+ \left[(I_{0^+}^{1-\gamma}\phi)(t)y(t)\right]_{t=0}^{t=T}.
 	\end{array}
 \end{equation}
 provided that all the terms involved in the computation make sense
\end{remark}
 Throughout this work without any mention, we will assume the following.

 \begin{assumption}\label{asump-beta-q}
 We assume that $\beta\in {C}[a,b]$,  $q\in L^\infty(a,b)$, and there are two constants $\beta_0>0$ and $q_0>0$ such that
 \begin{align*}
\beta(x)\ge \beta_0>0 \; \mbox{ for all } x\in [a,b],\;\;\mbox{ and }  q(x)\ge q_0>0 \;\mbox{ for a.e.  }\; x\in [a,b].
 \end{align*}
 \end{assumption}

We have the following result.

\begin{lemma}\label{lembil}
Let $0<\alpha\leq 1$ and
\begin{align}\label{esp-E}
D(\mathcal E):=\{u\in \Ha:\; I_{a^+}^{1-\alpha}(a)=0\}
\end{align}
be the closed subspace of $\Ha$.
For any $y,z\in D(\mathcal E)$. We define  $\mathcal E:D(\mathcal E) \times D(\mathcal E)\to\R$ by
\begin{align}\label{bil}
\mathcal E(y,z)=\int_a^b \beta(x)\mathbb{D}^\alpha_{a^+}y(x) \mathbb{D}^\alpha_{a^+}z(x)dx +\int_a^b q(x)y(x)z(x)dx.
\end{align}
Then, $\mathcal E$ is continuous, closed, and coercive.
\end{lemma}

As a direct consequence of Lemma \ref{lembil1} we have the following observation.

\begin{remark}\label{rem-213}
Let $0<\alpha\leq 1$.
\begin{enumerate}
\item Let $A$ be the self-adjoint operator in $L^2(a,b)$ associated with the form $\mathcal E$ in the sense that
\begin{equation}\label{op-form}
\begin{cases}
D(A)=\Big\{u\in D(\mathcal E):\;\exists\; f\in L^2(a,b):\; \mathcal E(u,v)=(f,v)_{L^2(a,b)}\;\forall\; v\in D(\mathcal E)\Big\}\\
A u=f.
\end{cases}
\end{equation}
More precisely, using the integration by parts formula \eqref{integ1&2} we obtain that
\begin{equation}\label{op-op}
\begin{cases}
D(A)=\Big\{u\in D(\mathcal E):\;\; \Dc^\alpha_{b^-}(\beta\Dr_{a^+}^\alpha u)\in L^2(a,b):\;   (\beta\Dr_{a^+}^{\alpha}u)(b)=0\Big\}\\
A u=\Dc^\alpha_{b^-}(\beta\Dr_{a^+}^\alpha u) +qu.
\end{cases}
\end{equation}
\item Let $D(\mathcal E)^\star$ be the dual of $D(\mathcal E)$ with respect to the pivot space $L^2(a,b)$ so that we have the continuous embedding $D(\mathcal E)\hookrightarrow L^2(a,b)\hookrightarrow D(\mathcal E)^\star$. Then, the operator $A$ can also be defined as a bounded operator from $D(\mathcal E)$ into $D(\mathcal E)^\star$ by
$$\langle Au,v\rangle_{D(\mathcal E)^\star, D(\mathcal E)}=\mathcal E(u,v),\;\; u,v\in D(\mathcal E).$$

\item  If $1/2<\alpha<1$,  then it follows from the embedding  \eqref{a6} that the operator $A$ has a compact resolvent. Let $(\mu_n)_n$ be the eigenvalues of $A$ with associated eigenfunctions $(\phi_n)$. It follows from the coercivity and the nonnegativity of $\mathcal E$ that
\begin{align*}
0<\mu_1\leq\mu_2\leq \cdots\leq \mu_n\leq\cdots\;\mbox{ with } \lim_{n\to\infty}\mu_n=+\infty.
\end{align*}

\item We notice that we can ensure that $A$ has a compact resolvent only in the case $1/2<\alpha\leq 1$. We do not know if this is the case when $0<\alpha\leq 1/2$, but this is not relevant in the present paper.
\end{enumerate}
\end{remark}

Throughout the following, if $H$ is a Hilbert space we denote by $(\cdot,\cdot)_H$ the associated scalar product, and if $\mathbb X$ is a Banach space with dual $\mathbb X^\star$, we let $\langle\cdot,\cdot\rangle_{ \mathbb X^\star,\mathbb X}$ be their duality pairing.

\begin{lemma}\label{lemK} 
Let $\V$ be a Hilbert space and $\V^\star $ its dual. Then, the space
	\begin{equation}\label{defK}
	\mathbb K:=\left\{\varphi\in L^2((0,T);\V):\; \Dr_T^\gamma\varphi\in L^2((0,T);\V^\star)\right\}
	\end{equation}
	endowed with the norm
		\begin{equation}\label{normK}
		  \norm{\varphi}_{\mathbb K}=\left(\norm{\varphi}_{L^2((0,T);\V)}^2+\norm{\Dr_T^\gamma\varphi}_{L^2((0,T);\V^\star)}^2\right)^{1/2}
		\end{equation}
is a Hilbert space.
\end{lemma}

\begin{proof} Let $(\rho_n)$ be a Cauchy sequence in $\mathbb K$. Then  $(\rho_n)$  and $\left(\Dr_T^\gamma(\rho_n)\right)$ are respectively Cauchy sequence in the Hilbert spaces $L^2((0,T);\V)$ and $L^2((0,T);\V^\star)$, respectively.  Hence,
\begin{equation}\label{inter12}
\begin{array}{rllllll}
\rho_n&\to & \rho& \hbox{ strongly in }L^2((0,T);\V),\\
\Dr_T^\gamma(\rho_n)&\to & \zeta& \hbox{ strongly in }L^2((0,T);\V^\star).
\end{array}
\end{equation}
To conclude, we need to prove that $\Dr_T^\gamma(\rho)=\zeta $ in $L^2((0,T);\V^\star)$. Indeed,
  let $\phi$ be a test function. Then using the integration by parts formula \eqref{integ1} we have that,
 $$\begin{array}{rllllll}
 \dis \int_0^T\langle\Dr_T^\gamma(\rho_n)(t),\phi(t)\rangle_{\V^\star,\V}dt &=&
 \dis  \int_0^T\langle\rho_n(t),\mathcal{D}_t^\gamma\phi(t)\rangle_{\V,\V^\star}dt.
 \end{array}
 $$
 Passing to the limit as$n\to \infty$ in this latter identity while using \eqref{inter12}, we obtain that
 $$
 \dis \int_0^T\langle \zeta(t),\phi(t)\rangle_{\V^\star,\V}dt=
 \dis  \int_0^T\langle\rho(t),\mathcal{D}_t^\gamma\phi(t)\rangle_{\V,\V^\star}dt.
 $$
 Using again the  integration by parts  formula \eqref{integ1} yields
 $$\begin{array}{rllllll}
 \dis \int_0^T\langle \zeta(t),\phi(t)\rangle_{\V^\star,\V}dt=
\dis \int_0^T\langle\Dr_T^\gamma(\rho)(t),\phi(t)\rangle_{\V^\star,\V}dt.
 \end{array}
 $$
 Thus, $\Dr_T^\gamma(\rho)=\zeta $ in $L^2((0,T);\V^\star)$ and the proof is finished.
\end{proof}

\begin{remark}\label{rem28} 
We observe the following,
\begin{enumerate}
\item Note that if $\rho\in \mathbb{K}$, then from  Remark \ref{rem27}, we obtain that $I^{1-\gamma}_{T}(\rho)\in L^2((0,T);\V)$. Therefore, using the fact that $\Dr_T^{\gamma}\phi=\dis -\frac{\partial }{\partial t }I^{1-\gamma}_{T}\phi\in L^2((0,T);\V^\star)$, we can deduce that $I^{1-\gamma}_{T}(\rho)\in \mathbb{W}((0,T);\V),$ where
\begin{equation}\label{defW0T}
	\mathbb{W}((0,T);\V):= \left\{\zeta \in L^2(0,T;\mathbb{V}): \dis\frac{\partial \zeta}{\partial t} \in L^2\left((0,T);\mathbb{V}^\star\right)\right\}.
	\end{equation}
	\item If   if $\mathbb{Y}$ is a Hilbert space that can be identified with its dual $\mathbb{Y}^\star$ and we have the continuous embeddings $\mathbb{V}\hookrightarrow \mathbb{Y}=\mathbb Y^\star\hookrightarrow \mathbb{V}^\star,$
	then using \cite{lions1971}
  we have the continuous embedding
  $$W(0,T;\mathbb{V})\hookrightarrow C([0,T];\mathbb{Y}).$$
This implies that
	\begin{equation}\label{contWTA}
	I^{1-\gamma}_{T}(\rho)\in  C([0,T];\mathbb{Y}).
	\end{equation}
\end{enumerate}
\end{remark}

\section{The state equation in a single edge} \label{oneedge}
In this section, we are  concerned with the following Sturm-Liouville problem  in a single edge:

\begin{equation}\label{pa1}
\left\{\begin{array}{lllllllll}
\mathcal{D}_t^\gamma y+\dis  \Dc_{b^-}^\alpha\,(\beta\, \Dr_{a^+}^\alpha y)+q\,y&=&f&\hbox{ in }& Q:={(0,T)\times (a,b)},\\
\dis (I_{a^+}^{1-\alpha} y)(\cdot,a^+)&=& 0&\hbox{ in }& (0,T),\\
\dis (\beta\Dr_{a^+}^{\alpha}y)(\cdot,b^-)&=&v&\hbox{ in }& (0,T),\\
y(0,\cdot)&=&y^0&\hbox{ in } &(a,b),
\end{array}
\right.
\end{equation}
where $f\in L^2(Q)$, $v\in L^2(0,T)$ and $y^0\in L^2(a,b)$.

Throughout the following, without any mention $D(\mathcal E)$ is the space introduced in \eqref{esp-E} and $\mathcal E$ is the bilinear form given in \eqref{bil}.


\subsection{The homogeneous case in a single edge}

In this section, we are concerned with the existence and regularity results of homogeneous fractional Sturm–Liouville parabolic equations of the following type:
\begin{equation}\label{pa}
\left\{\begin{array}{rlllllllll}
\mathcal{D}_t^\gamma y+\dis  \Dc_{b^-}^\alpha\,(\beta\, \Dr_{a^+}^\alpha y)+q\,y&=&f&\hbox{ in }& Q,\\
\dis (I_{a^+}^{1-\alpha} y)(\cdot,a^+)&=& 0&\hbox{ in }& (0,T),\\
\dis (\beta\Dr_{a^+}^{\alpha}y)(\cdot,b^-)&=&0&\hbox{ in }& (0,T),\\
y(0,\cdot)&=&y^0&\hbox{ in } &(a,b).
\end{array}
\right.
\end{equation}
It follows from the characterization \eqref{op-op} of the operator $A$ that \eqref{pa} can be rewritten as the following abstract Cauchy problem:
\begin{equation}\label{ACP}
\mathcal{D}_t^\gamma y+ Ay= f\;\;\mbox{ in } Q,\qquad
y(0,\cdot)=y^0 \mbox{ in } (a,b).
\end{equation}


\begin{definition}
A function $y$ is said to be a weak solution of \eqref{pa} if the following assertions hold:
\begin{itemize}
\item Regularity:
$\dis y\in C([0,T];L^2(a,b))\cap L^2((0,T);D(\mathcal E))\;\mbox{ and }\; \Dc_t^\gamma y\in L^2((0,T);D(\mathcal E)^\star).$
\item Initial condition:
$\dis y(0,\cdot)=y^0 \;\mbox{ a.e. in }\;(a,b).$
\item Variational identity: for every $\varphi\in D(\mathcal E)$ and a.e. $t\in (0,T)$,  
$$\dis\langle \Dc_t^\gamma y, \varphi\rangle_{D(\mathcal E)^\star, D(\mathcal E)}+\mathcal E(y(t,\cdot),\varphi)=\langle f(t,\cdot),\varphi)_{D(\mathcal E)^\star, D(\mathcal E)}.$$

\end{itemize}
\end{definition}

We have the following  result.

\begin{theorem}\label{Theo1}
Let  $0<\gamma,\alpha\leq 1$. Then, for every $f\in L^2((0,T);D(\mathcal E)^\star)$ and $y^0\in L^2(a,b)$, the system \eqref{pa} (hence,  \eqref{ACP}) has a unique weak solution $y$ and there is a constant $C>0$ such that
\begin{equation}\label{esti1}
	\norm{y}_{{C}([0,T];L^2(a,b))} +\|y\|_{L^2((0,T);D(\mathcal E))}\leq C\left[\norm{y^0}_{L^2(a,b)}+\norm{f}_{L^2((0,T);D(\mathcal E)^\star)}\right].
\end{equation}
If $1/2<\alpha<1$, then $y$ has the representation
\begin{align}
y(t,\cdot)=&\sum_{n=1}^\infty (y^0,\phi_n)_{L^2(a,b)} E_{\gamma,1}(-\mu_nt^\gamma)\phi_n\notag\\
&+\sum_{n=1}^\infty\left(\int_0^t\langle f(\tau,\cdot),\phi_n\rangle_{D(\mathcal E)^\star, D(\mathcal E)}(t-\tau)^{\gamma-1}E_{\gamma,\gamma}(-\mu_n(t-\tau)^\gamma)\;d\tau\right)\phi_n,
\end{align}
where $E_{\gamma,1}$ and $E_{\gamma,\gamma}$ are the Mittag-Leffler functions given in  \eqref{eq23}, and $(\phi_n)$ denote the eigenfunctions of $A$ with associated eigenvalues $(\mu_n)$.
\end{theorem}

\begin{proof}
The proof is based on the so called subordination principle and the representation \eqref{ACP}.  We refer to \cite{Baz,Ga-Wa-2021,yama2,War-SIAM} and their references for more details.
\end{proof}

\subsection{The non-homogeneous case  in a single edge}

We first consider the following auxiliary evolution equation:

\begin{equation}\label{paaux}
\left\{\begin{array}{rlllllllll}
\mathbb{D}_t^\gamma p+\dis  \Dc_{b^-}^\alpha\,(\beta\, \Dr_{a^+}^\alpha p)+q\,p&=&g&\hbox{ in }& Q,\\
\dis (I_{a^+}^{1-\alpha} p)(\cdot,a^+)&=& 0&\hbox{ in }& (0,T),\\
\dis (\beta\Dr_{a^+}^{\alpha}p)(\cdot,b^-)&=&0&\hbox{ in }& (0,T),\\
I_t^{1-\gamma}p(0,\cdot)&=&p^0&\hbox{ in } &(a,b).
\end{array}
\right.
\end{equation}
It also follows from  \eqref{op-op} that  \eqref{paaux} can be rewritten as the following abstract Cauchy problem:
\begin{equation}\label{ACPaux}
\mathbb{D}_t^\gamma p+ Ap= g\;\;\mbox{ in } Q,\qquad
I_t^{1-\alpha}p(0,\cdot)=p^0 \mbox{ in } (a,b).
\end{equation}


\begin{definition}
Let $g\in L^2((0,T);D(\mathcal E)^\star)$ and $p^0\in L^2(a,b)$.  A function $p$  is said to be a weak solution of \eqref{paaux}  if the following assertions hold:
\begin{itemize}
\item Regularity:
$\dis p\in L^2((0,T);D(\mathcal E)),\;  \mathbb D_t^\gamma p\in L^2((0,T);D(\mathcal E)^\star)~\mbox{ and }~\; I_{t}^{1-\gamma} p\in C([0,T];L^2(a,b)).$
\item Initial condition:
$\dis I_t^{1-\gamma}p(0,\cdot)=p^0 \;\mbox{ a.e. in }\;(a,b).$
\item Variational identity: for every $\varphi\in D(\mathcal E)$ and a.e. $t\in (0,T)$,  
$$\dis \langle \Dr_t^\gamma p(t,\cdot), \varphi\rangle_{D(\mathcal E)^\star, D(\mathcal E)}+\mathcal E(p(t,\cdot),\varphi)=\langle g(t,\cdot),\varphi\rangle_{D(\mathcal E)^\star, D(\mathcal E)}.$$
\end{itemize}
\end{definition}

We have the following existence result.

\begin{theorem}\label{Theo1aux}
Let  $0<\gamma,\alpha\leq 1$. Then,  for every $g\in L^2((0,T);D(\mathcal E)^\star)$ and $p^0\in L^2(a,b)$, the system \eqref{paaux} (hence,  \eqref{ACPaux}) has a unique weak solution $p$ and there is a constant $C> 0$ such that
\begin{align}\label{estim_p}
&	\norm{p}_{L^2((0,T);D(\mathcal E))}+ \norm{ \mathbb D_t^\gamma p}_{L^2((0,T);D(\mathcal E)^\star)}+\norm{I^{1-\gamma}_{t}p}_{C([0,T];L^2(a,b))}\notag\\
	\leq &C \left(\norm{g}_{L^2(0,T;D(\mathcal{E})^{\star})}+\norm{p^0}_{L^2(a,b)}\right).
\end{align}
If $1/2<\alpha\leq 1$, then $p$ has the representation
\begin{align*}
p(t,\cdot)=&\sum_{n=1}^\infty (p^0,\phi_n)_{L^2(a,b)} E_{\gamma,\gamma}(-\mu_nt^\gamma)\phi_n\notag\\
&+\sum_{n=1}^\infty\left(\int_0^t\langle g(\tau,\cdot),\phi_n\rangle_{D(\mathcal E)^\star, D(\mathcal E)}(t-\tau)^{\gamma-1}E_{\gamma,\gamma}(-\mu_n(t-\tau)^\gamma)\;d\tau\right)\phi_n.
\end{align*}
\end{theorem}

\begin{proof}
Since an abstract version of this result is contained in \cite{Baz,War-SIAM} and the references therein, we just give the main ideas of the proof.   Also we only consider the case $1/2<\alpha\leq 1$. The case $0<\alpha\leq 1/2$ follows with a similar abstract representation.

Let $(\phi_n)$ be the orthonormal basis defined in \eqref{op-op}.
Taking the duality map between  \eqref{ACPaux} and $\phi_n$,
setting $p_n(t):=(p(t), \phi_n)_{L^2(a,b)},\, g_n(t):=\langle g(t),\phi_n\rangle_{D(\mathcal{E})^\star, D(\mathcal{E})}$ and $p^0_n:=(p^0,\phi_n)_{L^2(a,b)}$,  we obtain that $p_n$ is a solution of the following  fractional ODE:
\begin{equation}\label{ACPaux2}
	\left\{\begin{array}{lllll}
		\mathbb{D}_t^{\gamma}p_n(t)+\mu_np_n(t)&=g_n(t)~~&\text{for }~~t\in (0,T),\\ I^{1-\gamma}_tp_n(0)&=p_n^0.~~&
	\end{array}\right.
\end{equation}
Using Laplace transform we get that $\eqref{ACPaux2}$ has a unique solution whose Laplace transform is given by
\begin{equation*}
	\hat{p}_n(s)=\frac{p^0_n}{s^{\gamma}+\mu_n}+\frac{\hat{g}_n(s)}{s^{\gamma}+\mu_n}.
\end{equation*}
Inverting the Laplace transform,  we obtain that
\begin{equation}\label{aux3}
	p_n(t)=t^{\gamma-1}E_{\gamma,\gamma}(-\mu_nt^{\gamma})p_n^0+\int_{0}^t(t-s)^{\gamma-1}E_{\gamma,\gamma}(-\mu_n(t-s)^{\gamma})g_n(s)\,\mathrm{d}s.
\end{equation}
Let
\begin{equation}
	p_m(t,\cdot)=\sum_{n=1}^{m}\left\{t^{\gamma-1}E_{\gamma,\gamma}(-\mu_nt^{\gamma})p_n^0+\int_{0}^t(t-s)^{\gamma-1}E_{\gamma,\gamma}(-\mu_n(t-s)^{\gamma})g_n(s)\,\mathrm{d}s\right\}\varphi_n
\end{equation}
and
$$
	p(t,\cdot)=\sum_{n=1}^{\infty}p_n(t)\phi_n=\sum_{n=1}^{\infty}\left\{t^{\gamma-1}E_{\gamma,\gamma}(-\mu_nt^{\gamma})p_n^0+\int_{0}^t(t-s)^{\gamma-1}E_{\gamma,\gamma}(-mu_n(t-s)^{\gamma})f_n(s)\,\mathrm{d}s\right\}\phi_n.$$
Using classical results on series, one shows that,  as $m\to\infty$,
\begin{equation*}
	\begin{array}{ll}
\dis	p_m&\to p \,~~\text{strongly in}\,~~ L^2((0,T);D(\mathcal E))\\
\mathbb D_t^\gamma p_m&\to p \,~~\text{weakly in}\,~~ L^2((0,T);D(\mathcal E)^\star)\\
\dis I^{1-\gamma}_t p_m&\to I^{1-\gamma}_t p \,~~\text{strongly in}\,~~C([0,T];L^2(a,b))\\
\mathcal E(p_m,\varphi)&\to \mathcal E(p,\varphi) \;\mbox{ for every }\varphi\in D(\mathcal E).
	\end{array}
\end{equation*}
Combining all these facts we obtain that $p$ also satisfies the variational equality. In addition, as $m\to\infty$,
 $$\dis	I^{1-\gamma}_{t}p_m(0,\cdot)=\sum_{i=1}^m p_i^0\varphi_i\to \sum_{i=1}^{\infty} p_i^0\varphi_i=p^0.$$
 We can deduce that $p$ is the unique  weak  solution of \eqref{paaux}.  The estimate \eqref{estim_p} is easy to establish by using the estimate \eqref{Mitl3}. The proof is finished.
\end{proof}

Now, we consider the non-homogeneous problem given in \eqref{pa1}.

\begin{definition}\label{defweak}
	Let $0<\gamma,\alpha \leq 1$, $f\in L^2(Q)$,   $y^0\in L^2(\Omega$, and  $v\in L^2(0,T)$. A function $y\in L^2(Q)$ is said to be a very-weak solution of \eqref{pa1} (or a solution by transposition) if the following identity holds:
	\begin{equation*} 
		\begin{array}{ll}&\dis\int_{Q}y\left(\D_T^{\gamma}\phi+\Dc_{b^-}^{\alpha}(\beta\Da\phi)+q\phi\right)\,dxdt\\&=\dis\int_{Q}f\phi\,dxdt+\int_{a}^by^0I^{1-\gamma}_T\phi(0,x)\mathrm{d}x+\int_{0}^Tv(t) I^{1-\alpha}_{a^+}\phi(t,b^-)\mathrm{d}t,
		\end{array}
	\end{equation*}
	 for all $\phi\in\Psi$, where
\begin{equation}\label{defPhi}
	\Psi:=\left\{\phi\in L^2((0,T);D(A)):\;\D_T^{\gamma}\phi\in L^2(Q),\; I^{1-\gamma}_{a^+}\phi\in C([0,T];L^2(a,b)),\; I^{1-\gamma}_{a^+}\phi(T,\cdot)=0\right\}.
\end{equation}
\end{definition}

We have the following result.

\begin{theorem} Let $0<\gamma,\alpha\leq 1$, $f\in L^2(Q)$, $y^0\in L^2(\Omega)$, and $v\in L^2(0,T)$. Then,  there exists a unique very-weak solution {$y\in L^2(Q)$} of \eqref{pa1} in sense of Definition \ref{defweak}. Moreover,  there is a constant $C>0$ such that
	\begin{equation*}
		\norm{y}_{L^2(Q)}\leq C \left(\norm{f}_{L^2(Q)}+\norm{y^0}_{L^2(\Omega)}+\norm{v}_{L^2(0,T)}\right).
	\end{equation*}
\end{theorem}

This result will be proved for the general case of the system \eqref{pa7} in Section \ref{sec-NHGG}.

\section{The state equation on general star graphs}\label{graph}

Here,  we are concerned with the existence and regularity of solutions to the state equation \eqref{pa7}.
It is clear that the first $m-1$ controls $v^i$ are the Dirichlet controls, while the $v^i , i = m + 1, \cdots , N$,  are the Neumann controls.

We set
$$
\mathbb H_{a}^\alpha:=\dis \prod_{i= 1}^N H_{a^+}^\alpha (a,b_i).$$
We endow $\mathbb L^2$ and $\mathbb{H}_{a}^\alpha$  with the norms given, respectively, by
\begin{equation*}
	\|\rho\|^2_{\mathbb L^2}=\dis \sum_{i=1}^N\|\rho^i\|^2_{L^2(a,b_i)},\quad \rho=(\rho^i)_i\in \mathbb L^2,
\end{equation*}
and
\begin{equation*} 
	\|\rho\|^2_{\mathbb{H}_{a}^\alpha}=\dis \sum_{i=1}^N\|\rho^i\|^2_{H^\alpha_{a^+}(a,b_i)}=\sum_{i=1}^N\left(
	\norm{\rho^i}^2_{L^2(a,b_i)}+\norm{\D^{\alpha}_{a^+}\rho^i}^2_{L^2(a,b_i)}\right),\quad \rho=(\rho^i)_i\in \mathbb{H}_{a}^\alpha.
	\end{equation*}

\begin{remark}\label{rem-cp}
It follows from \eqref{a6} that if $1/2<\alpha\leq 1$, then the embedding $\mathbb H_{a}^\alpha\hookrightarrow \mathbb L^2$ is compact.
\end{remark}

Throughout the following without any mention we assume that the coefficients satisfy the following.

\begin{assumption}\label{asump-beta-q1}
 We assume that for $i=1,\ldots, N$, $\beta^i\in C([a,b_i])$,  $q^i\in L^\infty(a,b_i)$, and there are constants $\beta_0^i>0$ and $q_0^i>0$ such that
 \begin{align*}
\beta^i(x)\ge \beta_0^i>0 \; \mbox{ for all } x\in [a,b_i],\;\;\mbox{ and }  q^i(x)\ge q_0^i>0 \;\mbox{ for a.e.  }\; x\in [a,b_i].
 \end{align*}
 \end{assumption}
	From now on, we set
\begin{subequations}\label{notations}
	\begin{alignat}{11}
		\underline{q^0}&:=&\dis \min_{1\leq i\leq n}q^{i},\quad
		\underline{\beta^0}&:=&\dis\min_{1\leq i\leq n}\beta^{i},\label{notations1}\\
		\overline{q}&:=&\dis \max_{1\leq i\leq n}\|q^i\|_{\infty},\quad
		\overline{\beta}&:=&\dis\max_{1\leq i\leq n}\|\beta^i\|_{\infty},\label{notations2}
	\end{alignat}
\end{subequations}
where $\|\beta^i\|_{\infty}:=\dis \max_{x\in[a,b_i]}|\beta^i(x)| \hbox{ and } \dis \|q^i\|_{\infty}:= \sup_{x\in(a,b_i)}|q^i(x)|.$

As for the case of a single edge studied in the previous section, we need some existence and regularity results for the homogeneous case.

\subsection{The homogeneous case in a general star graph}	

We  consider the following fractional Sturm-Liouville boundary value problem on a general star graph:
\begin{equation}\label{pa7ad}
	\left\{
	\begin{array}{lllllllllllllllll}
		\displaystyle	\Dc_t^{\gamma}\rho^i+\mathcal{D}_{b_i^-}^\alpha(\beta^i\mathbb{D}_{a^+}^\alpha \rho^i)+q^i\rho^i&=&f^i&\text{in}& Q_i,\,i=1,\dots, N,\\
		\displaystyle	I_{a^+}^{1-\alpha}\rho^i(\cdot,a^+)-I_{\cdot,a^+}^{1-\alpha}\rho^j(\cdot,a^+)&=&0&\text{in}&(0,T),~i\neq j=1,\dots, N,\\
		\displaystyle	\sum_{i=1}^N\beta^i(a)\mathbb{D}_{a^+}^\alpha \rho^i(\cdot,a^+)&=&0&\text{in}& (0,T), \\
		\displaystyle	I_{a^+}^{1-\alpha} \rho^i(\cdot,b_i^-)&=&0&\text{in}& (0,T), \; i=1,\dots, m\\
		\displaystyle	\beta^i(b_i)\mathbb{D}_{a^+}^\alpha \rho^i(\cdot,b_i^-)&=&0&\text{in}&(0,T), \; i=m+1,\dots,N,\\
		\displaystyle	\rho^i(0,\cdot) &=&\rho^{0,i}&\text{in}& (a,b_i),~~i=1,\dots,N,
	\end{array}
	\right.
\end{equation}
where $f=(f^i)_i\in L^2((0,T); \mathbb L^2)$ and $\rho^0=(\rho^{0,i})_i\in \mathbb L^2$.

\begin{remark}\label{lembil1}
	Let $0<\alpha\leq 1$ and let
\begin{align}\label{DefVstar}
\V:=\dis \Big\{\rho:=(\rho^i)_i\in\mathbb H_a^\alpha: \;  (I^{1-\alpha}_{a^+}\rho^i)(a)-(I^{1-\alpha}_{a^+}\rho^j)(a)=0, i\neq j,\,\,i,j=1,\ldots, N, \notag\\
\qquad \qquad\qquad\mbox{ and }\; (I_{a^+}^{1-\alpha} \rho^i)(b_i^-)= 0, \; i=1,\dots,m\Big\},
\end{align}
be endowed with the norm given by
\begin{equation}\label{normstar}
\| \rho\|^2_{\V}:=\sum_{i=1}^N\|\rho^i\|^2_{H^\alpha_{a^+}(a,b_i)}=\|\rho\|_{\mathbb H_a^\alpha}^2
\end{equation}
which is a closed subspace of  $\mathbb{H}_{a}^{\alpha}$. We let $\V^\star$ denote the dual of $\V$ with respect to the pivot space $\mathbb L^2$ so that we have the continuous and dense embeddings: $\V\hookrightarrow\mathbb L^2\hookrightarrow \V^\star$.
We define $\mathbb F:\V \times \V\to\R$ by
	\begin{align*}\label{bil1}
		\mathbb{F}(\rho,\phi):=&
		\dis \sum_{i=1}^N\int_{Q_i} (\beta^i(x)\Dr_{a^+}^\alpha \phi^i )(x,t) \Dr_{a^+}^\alpha (\rho^i))(x,t)  dx dt\\
		&+ \dis \sum_{i=1}^N\int_{Q_i}q^i(x)\phi^i(x,t)\rho^i(x,t) dx\,dt.
	\end{align*}
\begin{enumerate}
\item Then, $\mathbb F$ is continuous, closed and coercive.  Let $\mathcal A$ be the selfadjoint operator  in $\mathbb L^2$ associated  with $\mathbb F$ in the sense that
\begin{equation}
\begin{cases}
D(\mathcal A)=\{U\in \V:\; \exists F\in \mathbb L^2:\; \mathbb F(U,V)=(F,V)_{\mathbb L^2}\;\;\forall\; V\in \V\}\\
\mathcal AU=F.
\end{cases}
\end{equation}
For $\rho=(\rho^i)_i\in D(\mathcal A)$, we let
\begin{equation}
	\mathcal{A}\rho= ((\mathcal{A}\rho)^i)_i,\;\; i=1,\ldots, N.
\end{equation}
Each component is given by
\begin{equation}\label{defA}
	(\mathcal{A}\rho)^i=\Dc^{\alpha}_{b^-}(\beta^i\Da \rho^i)+q^i\rho^i,\;\;\;~~~~~~~ i=1,\cdots,N.
\end{equation}
\item If $1/2<\alpha\leq 1$, then it follows from Remark \ref{rem-cp} that the operator $\mathcal A$ has a compact resolvent. Let $(\lambda_k)_k$ be the eigenvalues of $\mathcal A$ with associated eigenfunctions $(\varphi_k)=(\varphi_{k}^i)_{1\leq i\leq N}$. It follows from the coercivity and the nonnegativity of $\mathbb F$ that
\begin{align*}
			0<\lambda_1\leq\lambda_2\leq \cdots\leq \lambda_k\leq\cdots\;\mbox{ with } \lim_{k\to\infty}\lambda_k=+\infty.
		\end{align*}
\end{enumerate}
		\end{remark}
		
Using Remark \ref{lembil1},  the system \eqref{pa7ad} can be rewritten as the  abstract Cauchy problem
		\begin{equation}\label{ACP1}
				\Dc_t^\gamma \rho+ \mathcal{A}\rho= f\quad \text {in  }~(0,T), \;\;\;\;\rho(0)=\rho^0.
		\end{equation}
		
	\begin{definition}
		A function $\rho$ is said to be a weak solution of \eqref{pa7ad} if the following assertions hold:
		\begin{itemize}
			\item Regularity:
			$\dis \rho\in C([0,T];\mathbb{L}^2)\cap L^2((0,T);\V)\;\mbox{ and }\; \Dc_t^\gamma \rho\in L^2((0,T);\V^\star).$
			\item Initial condition:
			$\dis  \rho^i(0,\cdot)=\rho^{i,0 }\;\mbox{ in }\;(a,b_i),~i=1,\cdots, N$.
			\item Variational identity: 	for every $\varphi\in \V$ and a.e. $t\in (0,T)$,  
			$$
			\dis\langle \Dc_t^\gamma \rho, \varphi\rangle_{\V^\star, \V}+\mathbb F(\rho(t,\cdot),\varphi)
			=\langle f(t,\cdot),\varphi)_{\V^\star\V}.$$
	\end{itemize}
	\end{definition}
	
	We have the following existence result which is also a direct consequence of the abstract results contained in \cite{Baz, Ga-Wa-2021}, the subordination principle,   and the representation \eqref{ACP1}.
	
	\begin{theorem}\label{Theo12}
		Let  $0<\gamma,\alpha\leq 1$. Then, for every $f\in L^2((0,T);\mathbb{V}^\star)$ and $\rho^0\in \mathbb{L}^2$, the system \eqref{pa7ad} (hence, the Cauchy problem \eqref{ACP1} has a unique weak solution $\rho$ and there is a constant $C>0$ such that
\begin{equation*}
	\norm{\rho}_{{C}([0,T];\mathbb L^2)} +\|\rho\|_{L^2((0,T);\V)}\leq C\left[\norm{y^0}_{\mathbb L^2}+\norm{f}_{L^2((0,T);\V^\star)}\right].
\end{equation*}
If  $1/2<\alpha\leq 1$,  then $\rho$ has the representation
		\begin{align}
			\rho^i(t,\cdot)=&\sum_{k=1}^\infty (\rho^{i,0},\varphi_k^i)_{L^2(a,b_i)} E_{\gamma,1}(-\lambda_kt^\gamma)\varphi_k^i\notag\\
			&+\sum_{k=1}^\infty\left(\int_0^t\langle f^i(\tau,\cdot),\varphi_k^i\rangle_{\V_i^\star,\V_i}(t-\tau)^{\gamma-1}E_{\gamma,\gamma}(-\lambda_k(t-\tau)^\gamma)\;d\tau\right)\varphi_k^i.
		\end{align}
	\end{theorem}

\subsection{The non-Homogeneous case in a general star graph}\label{sec-NHGG}

We first consider the following dual problem:
\begin{equation}\label{paaux1}
	\left\{
	\begin{array}{lllllllllllllllll}
		\displaystyle	\Dr_t^{\gamma}y^i+\mathcal{D}_{b_i^-}^\alpha(\beta^i\mathbb{D}_{a^+}^\alpha y^i)+q^iy^i&=&g^i&\text{in}& Q_i,\,i=1,\dots, N,\\
		\displaystyle	I_{a^+}^{1-\alpha}y^i(\cdot,a^+)-I_{a^+}^{1-\alpha}y^j(\cdot,a^+)&=&0&\text{in}&(0,T),~i\neq j=1,\dots, N,\\
		\displaystyle	\sum_{i=1}^N\beta^i(a)\mathbb{D}_{a^+}^\alpha y^i(\cdot,a^+)&=&0&\text{in}& (0,T), \\
		\displaystyle	I_{a^+}^{1-\alpha} y^i(\cdot,b_i^-)&=&0&\text{in}& (0,T), \; i=1,\dots, m\\
		\displaystyle	\beta^i(b_i)\mathbb{D}_{a^+}^\alpha y^i(\cdot,b_i^-)&=&0&\text{in}&(0,T), \; i=m+1,\dots,N,\\
		\displaystyle	I^{1-\gamma}_ty^i(0,\cdot) &=&y^{0,i}&\text{in}& (a,b_i),~~i=1,\dots,N.
	\end{array}
	\right.
\end{equation}
The system \eqref{paaux1} can be also rewritten as the abstract Cauchy problem
\begin{equation}\label{ACPaux3}
		\mathbb{D}_t^\gamma y+ \mathcal Ay= g\;\,\hbox{ in } (0,T),\;\; I_t^{1-\alpha}y(0)=y^0.
\end{equation}

\begin{definition}\label{definition46}
	A function $y$ is said to be a weak solution of \eqref{paaux1}  if the following assertions hold:
	\begin{itemize}
		\item Regularity: $\; I_{t}^{1-\gamma} y\in C([0,T];\mathbb L^2),\; y\in L^2((0,T);\V)~\mbox{ and }~\Dr_t^\gamma y\in L^2((0,T);\V^\star).$
		\item Initial condition:
		$\dis I_t^{1-\lambda}y(0^+,\cdot)=y^0.$
		\item Variational identity: for every $\varphi\in \V$ and a.e. $t\in (0,T)$, 
		$$\dis \langle \Dr_t^\gamma y(t,\cdot), \varphi\rangle_{\V^\star, \V}+\mathbb F(y(t,\cdot),\varphi)=\langle g(t,\cdot),\varphi\rangle_{\V^\star, \V}.$$
	\end{itemize}
\end{definition}

We have the following existence result.

\begin{theorem}\label{Theo1aux1}
	Let  $0<\gamma,\alpha\leq 1$.
	For every $g\in L^2((0,T);\V^\star)$ and $y^0\in \mathbb L^2$, the system \eqref{paaux1} (hence,  \eqref{ACPaux3}) has a unique weak solution $y$.  Moreover, there is a constant $C\ge 0$ such that
\begin{equation}\label{estim_y1}
	\norm{y}_{L^2(0,T,\V)} +\norm{I^{1-\gamma}_{t}y}_{C([0,T];\mathbb L^2)}+	\norm{\Dr_t^{\gamma}y}_{L^2(0,T;\V^\star)}\leq C \left(\norm{g}_{L^2(0,T;\V^\star)}+\norm{y_0}_{\mathbb L^2}\right).
\end{equation}
If $1/2<\alpha\leq 1$, 	 then $y$ has the representation (where $i=1,2,\ldots,N$)
	\begin{align}
		y^i(t,\cdot)=&\sum_{k=1}^\infty (y^{0,i},\varphi_k^i)_{L^2(a,b_i)} E_{\gamma,\gamma}(-\lambda_kt^\gamma)\varphi_k^i\notag\\
		&+\sum_{k=1}^\infty\left(\int_0^t\langle g^i(\tau,\cdot),\varphi_k^i\rangle_{\V_i^\star, \V_i}(t-\tau)^{\gamma-1}E_{\gamma,\gamma}(-\lambda_k(t-\tau)^\gamma)\;d\tau\right)\varphi_k^i.
	\end{align}
\end{theorem}

\begin{proof}
The proof of this theorem follows the lines of the case of a single edge given in Theorem \ref{Theo1aux}. We omit the details.
\end{proof}

\begin{remark}\label{reg-sol-D}
We observe that in \eqref{paaux1} or  \eqref{ACPaux3} if $g\in L^2(0,T;\mathbb L^2)$ and $y^0\in D(\mathcal A)$, then the weak solution $y$ also enjoys the following additional regularity:
$$y\in C([0,T];D(\mathcal A)) \hbox{ and } \Dr_t^\gamma y^i\in L^2(Q_i)\;, i=1,\ldots, N.$$
\end{remark}

Next, for every $g=(g^i)_i\in L^2(0,T;\mathbb{L}^2)$ and $\phi^T=(\phi^{T,i})_i\in \mathbb{L}^2$, we consider the associated dual system given by
	\begin{equation}\label{pa21}
		\left\{
		\begin{array}{lllllllllllllllll}
			\displaystyle	\Dr_T^{\gamma}\phi^i+\mathcal{D}_{b_i^-}^\alpha(\beta^i\mathbb{D}_{a^+}^\alpha \phi^i)+q^i\phi^i&=&g^i&\text{in}& Q_i,\,i=1,\dots, N,\\
			\displaystyle	I_{a^+}^{1-\alpha}\phi^i(\cdot,a^+)-I_{a^+}^{1-\alpha}\phi^j(\cdot,a^+)&=&0&\text{in}&(0,T),~i\neq j=1,\dots, N,\\
			\displaystyle	\sum_{i=1}^N\beta^i(a)\mathbb{D}_{a^+}^\alpha \phi^i(\cdot,a^+)&=&0&\text{in}& (0,T), \\
			\displaystyle	I_{a^+}^{1-\alpha} \phi^i(\cdot,b_i^-)&=&0&\text{in}& (0,T), \; i=1,\dots, m\\
			\displaystyle	\beta^i(b_i)\mathbb{D}_{a^+}^\alpha \phi^i(\cdot,b_i^-)&=&0&\text{in}&(0,T), \; i=m+1,\dots,N,\\
			\displaystyle	I^{1-\gamma}_T\phi^i(T,\cdot) &=&\phi^{T,i}&\text{in}& (a,b_i),~~i=1,\dots,N.
		\end{array}
		\right.
	\end{equation}
	
\begin{definition} \label{definition47} A function $\phi$ is said to be a weak solution of \eqref{pa21}  if the following assertions hold:
	\begin{itemize}
		\item Regularity: $I_{T}^{1-\gamma} \phi\in C([0,T];\mathbb L^2),\; \phi\in L^2((0,T);\V)~\mbox{ and }~\Dr_T^\gamma \phi\in L^2((0,T);\V^\star).$
		\item Initial condition:
		$\dis I_T^{1-\lambda}\phi(T^-,\cdot)=\phi^T.$
		\item Variational identity: for every $\varphi\in \V$ and a.e. $t\in (0,T)$, 
		$$\dis \langle \Dr_T^\gamma \phi(t,\cdot), \varphi\rangle_{\V^\star, \V}+\mathbb F(\phi(t,\cdot),\varphi)=\langle g(t,\cdot),\varphi\rangle_{\V^\star, \V}.$$
	\end{itemize}
\end{definition}

\begin{corollary}\label{coroaux}
Let $0<\gamma,\alpha\leq1$,  $\phi^T=(\phi^{T,i})_i\in \mathbb{L}^2$ and $g=(g^i)_i\in L^2(0,T;\mathbb{L}^2)$. Then,  there exists a unique weak solution $\phi$ of \eqref{pa21}. Moreover,  there is a constant $C>0$ such that
	\begin{equation}\label{estaux}
		\norm{	\Dr_T^{\gamma}\phi}_{L^2((0,T);\V^\star)}+\norm{\phi}_{L^2(0,T;\V)}+\norm{I^{1-\gamma}_T\phi}_{C([0,T];\mathbb L^2)}\leq C\left(\norm{g}_{L^2(0,T;\mathbb L^2)}+\norm{\phi^T}_{\mathbb L^2}\right).
	\end{equation}
\end{corollary}

\begin{proof}
	Let  $\rho(x,T-t):=\phi(x,t)$. Then, for every $t\in (0,T)$ we have that
	\begin{equation*}
			I^{1-\gamma}_T \phi(t)=I^{1-\gamma}_t\rho(T-t)\,\, \hbox{ and }\; \D_T^{\gamma}\phi(t)=\D_t^{\gamma}\rho(T-t),
	\end{equation*}
	and $\rho$ satisfies the following system:
	\begin{equation}\label{paaux12}
	\left\{
	\begin{array}{lllllllllllllllll}
		\displaystyle	\Dr_t^{\gamma}\rho^i+\mathcal{D}_{b_i^-}^\alpha(\beta^i\mathbb{D}_{a^+}^\alpha \rho^i)+q^i\rho^i&=&g^i&\text{in}& Q_i,\,i=1,\dots, N,\\
		\displaystyle	I_{a^+}^{1-\alpha}\rho^i(\cdot,a^+)=I_{a^+}^{1-\alpha}\rho^j(\cdot,a^+)&=&0&\text{in}&(0,T),~i\neq j=1,\dots, N,\\
		\displaystyle	\sum_{i=1}^N\beta^i(a)\mathbb{D}_{a^+}^\alpha \rho^i(\cdot,a^+)&=&0&\text{in}& (0,T), \\
		\displaystyle	I_{a^+}^{1-\alpha} \rho^i(\cdot,b_i^-)&=&0&\text{in}& (0,T), \; i=1,\dots, m\\
		\displaystyle	\beta^i(b_i)\mathbb{D}_{a^+}^\alpha \rho^i(\cdot,b_i^-)&=&0&\text{in}&(0,T), \; i=m+1,\dots,N,\\
		\displaystyle	I^{1-\gamma}_t\rho^i(0,\cdot) &=&
\phi^{T,i}&\text{in}& (a,b_i),~~i=1,\dots,N.
	\end{array}
	\right.
\end{equation}
In view of Theorem \ref{Theo1aux1}, there exists a unique $\rho\in L^2((0,T);\V) $ solution of \eqref{paaux12}. Moreover,  from \eqref{estim_y1},  we have that there is a constant $C>0$ such that
	\begin{equation}\label{gf}
\norm{	\Dr_t^{\gamma}\rho}_{L^2((0,T);\V^\star)}+\norm{\rho}_{L^2(0,T,\V)}+\norm{I^{1-\gamma}_T\rho}_{C([0,T];\mathbb L^2)}\leq C\left(\norm{g}_{L^2(0,T;\mathbb L^2)}+\norm{\phi^T}_{\mathbb L^2}\right).
	\end{equation}
We can then deduce that \eqref{pa21} has a unique solution $\phi\in L^2((0,T);\V) $ satisfying the estimate \eqref{estaux}.
\end{proof}

Next, we consider the state equation \eqref{pa7}.

\begin{definition}\label{defweak1}
	Let $0<\gamma,\alpha\leq 1$,  {$f\in L^2(0,T;\mathbb{L}^2)$, $y^0\in \mathbb{L}^2$,  and $v\in\left( L^2(0,T)\right)^{N-1}$. A function $y\in L^2(0,T;\mathbb L^2)$} is said to be a very-weak solution of \eqref{pa7}  if the  identity
	\begin{equation}\label{weak1}
		\begin{array}{ll}
&\dis\sum_{i=1}^N\dis\int_{Q_i}y^i\left(\D_T^{\gamma}\phi^i+\Dc_{b_i^-}^{\alpha}(\beta^i\Da\phi^i)+q^i\phi^i\right)\,dxdt\\&=\dis \sum_{i=1}^N\int_{Q_i}f^i\phi^i\,dxdt+\sum_{i=1}^N\int_{a}^{b_i} y^{i,0}I^{1-\gamma}_T\phi^i(x,0)\mathrm{d}x\\&\dis-\sum_{i=2}^m\int_{0}^Tv^i\beta^i(b_i)\Da\phi^i(b_i^-)\,dt+\sum_{i=m+1}^N\int_{0}^Tv^i(t) I^{1-\alpha}_{a^+}\phi^i(b_i^-,t)\mathrm{d}t,
		\end{array}
	\end{equation}
holds for all $\phi=(\phi^i)_{i=1,\cdots,N}\in\Phi$, where
	\begin{equation}\label{defPhi1}
		\Phi=\left\{\phi\in L^2((0,T),D(\mathcal A)):\; \D_T^{\gamma}\phi^i\in L^2(Q_i),\; I^{1-\gamma}_{T}\phi\in C([0,T];\mathbb L^2) \mbox{ and } I^{1-\gamma}_{T}\phi(\cdot,T)=0\right\}.
	\end{equation}
\end{definition}

We have the following existence result.

\begin{theorem}\label{very-weak}
Let $0<\gamma,\alpha\leq1$,  {$f\in L^2(0,T;\mathbb L^2)$,}  $y^0\in \mathbb L^2$,  and $v\in \left(L^2(0,T)\right)^{N-1}$. Then,  there exists a unique very-weak solution $y\in L^2((0,T;\mathbb L^2)$ of \eqref{pa7} in the sense of Definition \ref{defweak1}. Moreover,  there is a constant $C>0$ such that
	\begin{equation}\label{estimation11}
		\norm{y}_{L^2(0,T;\mathbb L^2)}\leq C \left(\norm{f}_{L^2(0,T;\mathbb L^2)}+\norm{y^0}_{\mathbb L^2}+\norm{v}_{(L^2(0,T))^{N-1}}\right).
	\end{equation}
\end{theorem}

\begin{proof}
For every $g\in L^2(0,T;\mathbb{L}^2)$, we consider the associated dual system given by \eqref{pa21}.
In view of Corollary \ref{coroaux}, we have that $I_{T}^{1-\gamma} \phi\in C([0,T];\mathbb L^2),\; \phi\in L^2((0,T);\V)~\mbox{ and }~\Dr_T^\gamma \phi\in L^2((0,T);\V^\star).$ Moreover,  there is a constant $C>0$ such that \eqref{estaux} holds. In addition, using  Remark \ref{reg-sol-D}, and the fact that $I^{1-\gamma}_{T}\phi(\cdot,T)=0$ in $\Omega$, we have that $\phi\in\Phi$, where $\Phi$ is the space given in \eqref{defPhi1}.
The rest of the proof is done in several steps. \par
{\bf Step 1}: We claim that there is a constant $C_1>0$ such that
	\begin{equation*}
		\sum_{i=1}^N\norm{\beta^i(b^i)\Da\phi^i(\cdot,b_i^-)}^2_{L^2(0,T)}\leq C_1 \norm{g}^2_{L^2(0,T; \mathbb L^2)}.
	\end{equation*}
Indeed,  first of all observe that
\begin{equation*}
	I^\alpha_{a^+}\Dc_{b_i^-}^{\alpha}(\beta^i\Da\phi^i)(\cdot,x)=\beta^i(b^i)\Da\phi^i)(\cdot,b_i^-)-(\beta^i\Da\phi^i)(\cdot,x).
\end{equation*}
Applying $I^\alpha_{a^+}$ to the first equation in $\eqref{pa21}$, we have that
\begin{equation*}
I^\alpha_{a^+}\Dr_T^{\gamma}\phi^i+I^\alpha_{a^+}\mathcal{D}_{b_i^-}^\alpha(\beta^i\mathbb{D}_{a^+}^\alpha \phi^i)+I^\alpha_{a^+}q^i\phi^i=I^\alpha_{a^+}g^i.
		\end{equation*}
This implies that
\begin{equation*}
\beta^i(b^i)\Da\phi^i)(\cdot,b_i^-)=(\beta^i\Da\phi^i)-I^\alpha_{a^+}\Dr_T^{\gamma}\phi^i-I^\alpha_{a^+}q^i\phi^i+I^\alpha_{a^+}g^i.
\end{equation*}
Taking the norm and using Lemma \ref{l1} gives
\begin{align*}
		&\left(\sum_{i=1}^N\norm{\beta^i(b^i)\Da\phi^i(\cdot,b_i^-)}^2_{L^2(0,T)}\right)^{\frac1 2}\\&\leq \overline{\beta}\norm{\Da\phi}_{L^2(0,T;\mathbb L^2)}+\norm{I^\alpha_{a^+}\Dr_T^{\gamma}\phi}_{L^2(0,T;\mathbb L^2)}+ \dis \norm{I^\alpha_{a^+}q\phi}_{L^2(0,T;\mathbb L^2)}+  \norm{I^\alpha_{a^+}g}_{L^2(0,T; \mathbb L^2)}\\
		&\dis \leq  \overline{\beta}\norm{\Da\phi}_{L^2(0,T;\mathbb L^2)}+C\norm{\Dr_T^{\gamma}\phi}_{L^2(0,T;\mathbb L^2)}+ C\dis \norm{q\phi}_{L^2(0,T;\mathbb L^2)}+  C\norm{g}_{L^2(0,T; \mathbb L^2)}.
\end{align*}
Using the estimate \eqref{estim_y1},  we get that there is a constant $C>0$ such that
\begin{align*}
		&\left(\sum_{i=1}^N\norm{\beta^i(b^i)\Da\phi^i(\cdot,b_i^-)}^2_{L^2(0,T)}\right)^{\frac1 2}\\&\leq \dis \overline{\beta}\norm{\Da\phi}_{L^2(0,T;\mathbb L^2)}+C\norm{g}_{L^2(0,T;\mathbb L^2)}+ C\overline{q}\dis \norm{\phi}_{L^2(0,T;\mathbb L^2)}+  C\norm{g}_{L^2(0,T; \mathbb L^2)}\\&\dis\leq\left(\overline{\beta}^2+C^2\overline{q}^2\right)^{1/2}\norm{\phi}_{L^2(0,T,\V)}+C\norm{g}_{L^2(0,T; \mathbb L^2)}
		\\&\dis\leq\left(\overline{\beta}^2+C^2\overline{q}^2\right)^{1/2}\norm{g}_{L^2(0,T,\mathbb L^2)}+C\norm{g}_{L^2(0,T; \mathbb L^2)}\\&\leq C_1 \norm{g}_{L^2(0,T; \mathbb L^2)},
	\end{align*}
and we have shown the claim.

{\bf Step 3}: Next,  we endow $\Phi$ with the norm given by
	\begin{align}\label{norm-Phi}
		\norm{\phi}_{\Phi}^2=\sum_{i=1}^N\norm{\D_T^{\gamma}\phi^i}_{L^2(Q_i)}^2+\sum_{i=1}^N\norm{\beta^i(b^i)\Da\phi^i(\cdot,b_i^-)}^2_{L^2(0,T)}+\norm{\phi}_{L^2((0,T),D(\mathcal A)}^2.
	\end{align}
We consider the map
	$$\dis \mathcal{L}:\Phi\to L^2(0,T;\mathbb L^2),\,\;\phi\mapsto \mathcal{L}\phi=	\mathbb{D}_T^\gamma \phi+\dis  \mathcal A\phi.$$
It follows from Step 1 that $\mathcal{L}$ is an isomorphism.

We consider the linear functional $\mathcal{M}:\Phi\to \mathbb R$ given by
\begin{equation*}
	\begin{array}{ll}
		\mathcal{M}(\phi)=&\dis \sum_{i=1}^N\int_{Q_i}f^i\phi^i\,dxdt+\sum_{i=1}^N\int_{a}^{b_i} y^{i,0}I^{1-\gamma}_T\phi^i(0,x)\mathrm{d}x\\&\dis-\sum_{i=2}^m\int_{0}^Tv^i\beta^i(b_i)\Da\phi^i(t,b_i^-)\,dt+\sum_{i=m+1}^N\int_{0}^Tv^i(t) I^{1-\alpha}_{a^+}\phi^i(t,b_i^-)\mathrm{d}t.
	\end{array}.
\end{equation*}
Using Step 2, Lemmas \ref{l1},  \ref{trace} and Remark \ref{remtrace},  we have that there is a constant $C>0$ such that
\begin{align}\label{M0}
	&	\abs{\mathcal{M}(\phi)}\dis\leq \norm{f}_{L^2(0,T,\mathbb L^2)}\norm{\phi}_{L^2(0,T,\mathbb L^2)}+\norm{y^0}_{\mathbb L^2}\norm{I^{1-\gamma}_T\phi(0,\cdot)}_{\mathbb L^2}\notag\\
		&\dis+\norm{v}_{(L^2(0,T))^{N-1}}\left(\sum_{i=m+1}^N\norm{\Da\phi(\cdot,b^-_i)}_{L^2(0,T)}^2\right)^{\frac1 2}+\norm{v}_{(L^2(0,T))^{N-1}}\left(\sum_{i=m+1}^N\norm{I^{1-\alpha}_{a^+}\phi(\cdot,b^-_i)}_{L^2(0,T)}^2\right)^{\frac1 2}\notag\\
		&\leq \left(\norm{f}_{L^2(0,T,\mathbb L^2)}^2+\norm{y^0}_{\mathbb L^2}^2+\norm{v}_{(L^2(0,T))^{N-1}}^2\right)^{1/2}\notag\\
		&\times\left(\norm{\phi}_{L^2(0,T,\mathbb L^2)}^2+\norm{I^{1-\gamma}_T\phi(0,\cdot)}_{\mathbb L^2}^2+\dis\sum_{i=m+1}^N\norm{\Da\phi(\cdot,b^-_i)}_{L^2(0,T)}^2+\dis\sum_{i=m+1}^N\norm{I^{1-\alpha}_{a^+}\phi(\cdot,b^-_i)}_{L^2(0,T)}^2\right)^{1/2}\notag\\
		&\leq \left(\norm{f}_{L^2(0,T,\mathbb L^2)}^2+\norm{y^0}_{\mathbb L^2}^2+\norm{v}_{(L^2(0,T))^{N-1}}^2\right)^{1/2}\notag\\
		&\times\left((1+C)\norm{\phi}_{L^2(0,T,\V)}^2+\norm{I^{1-\gamma}_T\phi}_{L^2(0,T;\V)}^2+\dis\sum_{i=1}^N\norm{\beta^i(b^i)\Da\phi^i(\cdot,b_i^-)}^2_{L^2(0,T)}\right)^{1/2}\notag\\
		&\leq C\left(\norm{f}_{L^2(0,T,\mathbb L^2)}^2+\norm{y^0}_{\mathbb L^2}^2+\norm{v}_{(L^2(0,T))^{N-1}}^2\right)^{1/2}\norm{\phi}_{\Phi}.
\end{align}
We have shown that $\mathcal{M}$ is a continuous linear functional on $\Phi$.
Since $\mathcal{L}$ is an isomorphism, we can deduce  from \cite[Capter III Section 9.2]{lions1971} that there exists a unique $y\in L^2(0,T;\mathbb L^2)$ such that for all $\phi\in\Phi$,
\begin{equation*}
\sum_{i=1}^N\int_{Q_i}y^i\left(\mathbb{D}_T^\gamma \phi^i+ (\mathcal A\phi)^i\right)\;dxdt=\mathcal M(\phi).
\end{equation*}
That is,  $y$ is the unique very-weak solution of \eqref{pa7} in the sense of Definition \ref{defweak1}.

{\bf Step 4}: Now,  taking $y=g$ in \eqref{pa21}, using \eqref{weak1} and \eqref{M0}, we get
\begin{equation}\label{mw}
	\norm{y}^2_{L^2(0,T,\mathbb L^2)}\leq C\left(\norm{f}_{L^2(0,T;\mathbb L^2)}^2+\norm{y^0}_{\mathbb L^2}^2+\norm{v}_{(L^2(0,T))^{N-1}}^2\right)^{1/2}\norm{\phi}_{\Phi}.
\end{equation}
Combining  \eqref{estim_y1}-\eqref{mw}, we obtain that
\begin{equation*}
	\norm{y}^2_{L^2(0,T), \mathbb L^2}\leq C\left(\norm{f}_{L^2(0,T;\mathbb L^2)}^2+\norm{y^0}_{\mathbb L^2}^2+\norm{v}_{L^2(0,T)}^2\right)^{1/2}\norm{y}_{L^2(0,T;\mathbb L^2)},
\end{equation*}
and we have shown  \eqref{estimation11}.
The proof is finished.
\end{proof}

\section{Existence of minimizers and optimality conditions in a general star graph}\label{no-regret}

In this section we  are concerned with the existence of minimizers of the optimal control problem \eqref{op1}-\eqref{pa7} where recall that the functional $J$ is given by
{\begin{equation}\label{mp2}
	{J}(v,y^0)=\sum_{i=1}^{N}\int_{Q_i} \abs{y^i(v,y^0)-y_d^i}^2\,\mathrm{d}x\mathrm{d}t+\zeta\sum_{i=2}^{N}\int_0^T|v^i|^{2} \,\mathrm{d}t,
\end{equation}
where $y(v,y^0)=(y^i(v,y^0))_i$ is the weak solution of \eqref{pa7},  $y_d=(y_d^i)_i\in \mathbb{L}^2$, and $\mathcal U_{ad}$ is a closed and convex subset of $\left(L^2(0,T)\right)^{N-1}$.}

\begin{definition}\label{def-LR}
	We say that $v=(v^i)_i$ is a no-regret control of \eqref{op1}-\eqref{pa7} if $v$ is a solution of
	\begin{equation}\label{inte4}
		\inf_{v\in \left(L^2(0,T)\right)^{N-1} }\left[\sup_{y^0\in  \mathbb L^2} J(v,y^0)-J(0,y^0)\right].
	\end{equation}
\end{definition}

Next, we observe that  
$$y(v,y^0)=y(v,0)+y(0,y^0)-y(0,0),$$
where  $y(v,0)=(y^i(v,0))_i,\,y(0,y^0)=(y^i(0,y^0))_i$ and $y(0,0)=(y^i(0,0))_i$ satisfy respectively,   the following systems:
\begin{equation}\label{payV}
	\left\{
	\begin{array}{lllllllllllllllll}
		\displaystyle	\Dc_t^{\gamma}y^i(v,0)+\mathcal{D}_{b_i^-}^\alpha(\beta^i\mathbb{D}_{a^+}^\alpha y^i(v,0))+q^iy^i(v,0)&=&f^i&\text{in}& Q_i,\,i=1,\dots, N,\\
		\displaystyle	I_{a^+}^{1-\alpha}y^i(a^+;v,0)-I_{a^+}^{1-\alpha}y^j(a^+;v,0)&=&0&\text{in}&(0,T),~i\neq j=1,\dots, N,\\
		\displaystyle	\sum_{i=1}^N\beta^i(a)\mathbb{D}_{a^+}^\alpha y^i(a^+;v,0)&=&0&\text{in}& (0,T), \\
		\displaystyle	I_{a^+}^{1-\alpha} y^1(b_1^-;v,0)&=&0&\text{in}& (0,T),\\
		\displaystyle	I_{a^+}^{1-\alpha} y^i(b_i^-;v,0)&=&v^i&\text{in}& (0,T), \; i=2,\dots, m\\
		\displaystyle	\beta^i(b_i)\mathbb{D}_{a^+}^\alpha y^i(b_i^-;v,0)&=&v^i&\text{in}&(0,T), \; i=m+1,\dots,N,\\
		\displaystyle y^i(0,\cdot;v,0) &=&0&\text{in}& (a,b_i),~~i=1,\dots,N,
	\end{array}
	\right.
\end{equation}
and
\begin{equation}\label{payF}
	\left\{
\begin{array}{lllllllllllllllll}
	\displaystyle	\Dc_t^{\gamma}y^i(0,y^0)+\mathcal{D}_{b_i^-}^\alpha(\beta^i\mathbb{D}_{a^+}^\alpha y^i(0,y^0))+q^iy^i(0,y^0)&=&f^i&\text{in}& Q_i,\,i=1,\dots, N,\\
	\displaystyle	I_{a^+}^{1-\alpha}y^i(a^+;0,y^0)-I_{a^+}^{1-\alpha}y^j(a^+;0,y^0)&=&0 &\text{in}&(0,T),~i\neq j=1,\dots, N,\\
	\displaystyle	\sum_{i=1}^N\beta^i(a)\mathbb{D}_{a^+}^\alpha y^i(a^+;0,y^0)&=&0&\text{in}& (0,T), \\
	\displaystyle	I_{a^+}^{1-\alpha} y^i(b_i^-;0,y^0)&=&0&\text{in}& (0,T), \; i=1,\dots, m\\
	\displaystyle	\beta^i(b_i)\mathbb{D}_{a^+}^\alpha y^i(b_i^-;0,y^0)&=&0&\text{in}&(0,T), \; i=m+1,\dots,N,\\
	\displaystyle y^i(0,\cdot;0,y^0) &=&y^{0,i}&\text{in}& (a,b_i),~~i=1,\dots,N,
\end{array}
\right.
\end{equation}
and
\begin{equation}\label{payO}
	\left\{
\begin{array}{lllllllllllllllll}
	\displaystyle	\Dc_t^{\gamma}y^i(0,0)+\mathcal{D}_{b_i^-}^\alpha(\beta^i\mathbb{D}_{a^+}^\alpha y^i(0,0))+q^iy^i(0,0)&=&f^i&\text{in}& Q_i,\,i=1,\dots, N,\\
	\displaystyle	I_{a^+}^{1-\alpha}y^i(a^+;0,0)-I_{a^+}^{1-\alpha}y^j(a^+;0,0)&=&0&\text{in}&(0,T),~i\neq j=1,\dots, N,\\
	\displaystyle	\sum_{i=1}^N\beta^i(a)\mathbb{D}_{a^+}^\alpha y^i(a^+;0,0)&=&0&\text{in}& (0,T), \\
	\displaystyle	I_{a^+}^{1-\alpha} y^i(b_i^-;0,0)&=&0&\text{in}& (0,T), \; i=1,\dots, m\\
	\displaystyle	\beta^i(b_i)\mathbb{D}_{a^+}^\alpha y^i(b_i^-;0,0)&=&0&\text{in}&(0,T), \; i=m+1,\dots,N,\\
	\displaystyle y^i(0,\cdot;0,0) &=&0&\text{in}& (a,b_i),~~i=1,\dots,N.
\end{array}
\right.
\end{equation}
Since $f=(f^i)_i\in L^2((0,T); \mathbb L^2)$,  $y^0=(y^{0,i})_i\in \mathbb L^2$ and $v\in\left( L^2(0,T)\right)^{N-1}$, in view of Theorem \ref{very-weak} and Theorem \ref{Theo12} we have that $y(v,0)\in L^2((0,T;\mathbb L^2))$ and $\,y(0,y^0),  y(0,0)\in L^2((0,T);\V)\cap {C}([0,T];\mathbb L^2)$.

\begin{lemma}
For any $v\in(L^2(0,T))^{N-1}$ and $y^0\in  \mathbb L^2$, the functional $J$ satisfies
	\begin{equation}\label{Inter1}
		\begin{array}{lll}J(v,y^0)-J(0,y^0)&=&\dis J(v,0)-J(0,0)\\&&\dis + 2\langle y(v,0)-y(0,0) ,y(0,y^0)-y(0,0)\rangle_{L^2((0,T); \mathbb L^2)},
		\end{array}
	\end{equation} where $y(v,0),\,y(0,y^0)$ and $y(0,0)$ are solutions of \eqref{payV}, \eqref{payF} and \eqref{payO}, respectively.
\end{lemma}

\begin{proof}
Observing on the one hand that  $y(v,y^0)=y(v,0)+y(0,y^0)-y(0,0)$, and on the other hand that
$$
	{J}(v,0)=\dis \sum_{i=1}^{N}\int_{Q_i} \abs{y^i(v,y^0)-y_d^i}^2\,\mathrm{d}x\mathrm{d}t+\zeta\sum_{i=2}^{N}\int_0^T|v^i|^{2} \,\mathrm{d}t,
$$
$$
	{J}(0,y^0)=\sum_{i=1}^{N}\int_{Q_i} \abs{y^i(0,y^0)-y_d^i}^2\,\mathrm{d}x\mathrm{d}t,
$$
and
$$
	{J}(0,0)=\sum_{i=1}^{N}\int_{Q_i} \abs{y^i(0,0)-y_d^i}^2\,\mathrm{d}x\mathrm{d}t,$$
we can write
	\begin{align*}
	{J}(v,y^0)=&\dis\norm{y(v,y^0)-y_d}^2_{L^2((0,T); \mathbb L^2)}\,+\zeta\norm{v}^2_{(L^2(0,T))^{N-1}}\\
	=&\dis \norm{y(v,0)+y(0,y^0)-y(0,0)-y_d}^2_{L^2((0,T); \mathbb L^2)}\,+\zeta\norm{v}^2_{(L^2(0,T))^{N-1}}\\
	=&\dis \norm{y(v,0)-y(0,0)}^2_{L^2((0,T); \mathbb L^2)}+\norm{y(0,y^0)-y_d}^2_{L^2((0,T); \mathbb L^2)}\\
	&\dis +2\langle y(v,0)-y(0,0) ,y(0,y^0)-y_d\rangle_{L^2((0,T); \mathbb L^2)}+\zeta\norm{v}^2_{(L^2(0,T))^{N-1}}.
	\end{align*}
This implies that
\begin{align*}
		{J}(v,y^0)-J(0,y^0)=&\norm{y(v,0)-y(0,0)}^2_{L^2((0,T); \mathbb L^2)}\\
&+\dis 2\langle y(v,0)-y(0,0) ,y(0,y^0)-y_d\rangle_{L^2((0,T); \mathbb L^2)}+\zeta\norm{v}^2_{(L^2(0,T))^{N-1}}\\
=&\dis \norm{y(v,0)-y_d}^2_{ L^2((0,T); \mathbb L^2)}+\norm{y(0,0)-y_d}^2_{ L^2((0,T); \mathbb L^2)}\,\\
&-2\dis \langle y(v,0)-y_d ,y(0,0)-y_d\rangle_{L^2((0,T); \mathbb L^2)}\\
&+\dis
2\langle y(v,0)-y(0,0) ,y(0,y^0)-y_d\rangle_{L^2((0,T); \mathbb L^2)}+\zeta\norm{v}^2_{(L^2(0,T))^{N-1}}\\
=&\dis {J}(v,0)-J(0,0)\\
&-2\dis \langle y(v,0)-y(0,0) ,y(0,0)-y_d\rangle_{ L^2((0,T); \mathbb L^2)}\\
&+\dis
2\langle y(v,0)-y(0,0) ,y(0,y^0)-y_d\rangle_{ L^2((0,T); \mathbb L^2)}\\
=&\dis {J}(v,0)-J(0,0)\\
&+\dis
2\langle y(v,0)-y(0,0) ,y(0,y^0)-y(0,0)\rangle_{L^2((0,T); \mathbb L^2)}.
	\end{align*}
The proof is finished.
\end{proof}

We have the following result.

\begin{lemma}
For every $v\in( L^2(0,T))^{N-1}$ and  $y^0\in  \mathbb L^2$, we have that
	\begin{equation}\label{Inter2}
		{J}(v,f)-J(0,f)= J(v,0)-J(0,0)+2\sum_{i=1}^N\int_a^{b_i}	I^{1-\gamma}_T \phi^i(0,\cdot;v)y^{0,i}\,dx ,
	\end{equation}
	where $\phi(v):=(\phi^i(v))_i=\left(\phi^i(x,t;v)\right)_i$ is the weak solution of
	\begin{equation}\label{PaPhi}
	\left\{\begin{array}{lllllllll}
		\displaystyle	\Dr_T^{\gamma}\phi^i(v)+\mathcal{D}_{b_i^-}^\alpha(\beta^i\mathbb{D}_{a^+}^\alpha \phi^i(v)+q^i\phi^i(v)&=&y^i(v,0)-y^i(0,0)&\text{in}& Q_i,\,i=1,\dots, N,\\
		\displaystyle	I_{a^+}^{1-\alpha}\phi^i(a^+;v)-I_{a^+}^{1-\alpha}\phi^j(a^+;v)&=&0&\text{in}&(0,T),~i\neq j=1,\dots, N,\\
		\displaystyle	\sum_{i=1}^N\beta^i(a)\mathbb{D}_{a^+}^\alpha \phi^i(a^+;v)&=&0&\text{in}& (0,T), \\
		\displaystyle	I_{a^+}^{1-\alpha} \phi^i(b_i^-;v)&=&0&\text{in}& (0,T), \; i=1,\dots, m\\
		\displaystyle	\beta^i(b_i)\mathbb{D}_{a^+}^\alpha \phi^i(b_i^-;v)&=&0&\text{in}&(0,T), \; i=m+1,\dots,N,\\
		\displaystyle	I^{1-\gamma}_T \phi^i(T;v) &=&0&\text{in}& (a,b_i),~~i=1,\dots,N.
		\end{array}
	\right.
	\end{equation}
\end{lemma}

\begin{proof}
Since $(y^i(v,0)-y^i(0,0))\in L^2((0,T); L^2(a, b_i))$, it follows from  Corollary \ref{coroaux} that \eqref{PaPhi} has a unique weak solution $\phi$ satisfying $I_{T}^{1-\gamma} \phi\in C([0,T];\mathbb L^2),\; \phi\in L^2((0,T);\V)~\mbox{ and }~\Dr_T^\gamma \phi\in L^2((0,T);\V^\star).$    Moreover, there is a constant $C>0$ such that
	\begin{equation}\label{sam1}
		\norm{	\Dr_T^{\gamma}\phi}_{L^2((0,T);\V^\star)}+\norm{\phi }_{L^2(0,T;\V)}+\norm{I^{1-\gamma}_T\phi}_{C([0,T];\mathbb L^2)}\leq C\norm{y^{\tau}(v,0)-y(0,0)}_{L^2(0,T;\mathbb L^2)}.
	\end{equation}
Now, multiplying the first equation in \eqref{PaPhi} by $\left(y^i(0,y^0)-y^i(0,0)\right)$ and using the integration by parts formula given in Lemma \ref{lem0}, we obtain that
	\begin{align*}
			&\dis\sum_{i=1}^N \int_{Q_i} (y^i(v,0)-y^i(0,0))(y^i(0,y^0)-y^i(0,0))\,dxdt\\
			=&\dis \sum_{i=1}^N  \int_{Q_i}	\left(\mathbb{D}_T^\gamma \phi^i(v)+\dis  \Dc_{b_i^-}^\alpha\,(\beta_i\, \Dr_{a^+}^\alpha \phi^i(v))+q^i\,\phi^i(v)\right)\left((y^i(0,y^0)-y^i(0,0))\right)\,dxdt\\
			=&\dis \sum_{i=1}^N \int_{Q_i}	\left[\mathcal{D}_t^\gamma (y^i(0,y^0)-y^i(0,0))+\dis  \Dc_{b_i^-}^\alpha\,(\beta_i\, \Dr_{a^+}^\alpha (y^i(0,y^0)-y^i(0,0)))\right]\phi^i(v)\,dxdt\\
			\dis
&+\sum_{i=1}^N \int_{Q_i}q^i(y^i(0,y^0)-y^i(0,0))\phi^i(v)\,dxdt+
			\dis \sum_{i=1}^N \int_a^{b_i}	I^{1-\gamma}_T \phi^i(0,\cdot;v)y^{0,i}\,dx\\
			=&
			\dis \sum_{i=1}^N \int_a^{b_i}	I^{1-\gamma}_T \phi^i(0,\cdot;v)y^{0,i}\,dx.
		\end{align*}
Combining the latter identity with \eqref{Inter1}, we get \eqref{Inter2}. The proof is finished.
\end{proof}

\begin{remark}
We observe the following.
\begin{enumerate}
\item It follows from Definition \ref{def-LR} that $v$ is a no-regret control of \eqref{op1}-\eqref{pa7} if $v$ is a solution of
\begin{equation}\label{Op3}
	\inf_{v\in\left(L^2(0,T)\right)^{N-1} }\left[ J(v,0)-J(0,0)+\sup_{y^0\in  \mathbb L^2}\sum_{i=1}^{N}\int_a^{b_i}	I^{1-\gamma}_T \phi^i(0,\cdot;v)y^{0,i}\,dx \right].
\end{equation}
Since
	\begin{equation}\label{P1}
		\sup_{y^0\in  \mathbb L^2}\sum_{i=1}^{N}\int_a^{b_i}	I^{1-\gamma}_T \phi^i(0,\cdot;v)y^{0,i}\,dx =0
	\end{equation}
	or
	\begin{equation}\label{P2}
		\sup_{y^0\in  \mathbb L^2}\sum_{i=1}^{N}\int_a^{b_i}	I^{1-\gamma}_T \phi^i(0,\cdot;v)y^{0,i}\,dx =+\infty,
	\end{equation}
it follows that a no-regret control $v$ will exist only in the case \eqref{P1}. This means that a no-regret control belongs to the set
\begin{equation}\label{defU}\mathcal{U}=\left\{v\in \left(L^2(0,T)\right)^{N-1}:\;\, I^{1-\gamma}_T \phi^i(0,\cdot;v)\in \left( \mathbb L^2\right)^{\bot} \right\}.\end{equation}
\item The main difficulty arising in the characterization of the no-regret control  is that we do not know the structure of the set $\mathcal{U}$.  Therefore,   we relax the functional in \eqref{inte4}.  More precisely, for every $\tau>0$, we consider the following minimization problem:
\begin{equation}\label{Op4}
	\inf_{v\in\left(L^2(0,T)\right)^{N-1} }\left[\sup_{y^0\in \mathbb L^2}\left( J(v,y^0)-J(0,y^0)-\tau\norm{y^0}^2_{ \mathbb L^2}\right) \right],
\end{equation}
which in view of \eqref{Inter2} is equivalent to
	\begin{equation}\label{Op5}
	\inf_{v\in \left(L^2(0,T)\right)^{N-1} }\left[ J(v,0)-J(0,0)+2\sup_{y^0\in  \mathbb L^2}\left(\langle y^0, I^{1-\gamma}_T \phi^i(0,\cdot;v)\rangle_{ \mathbb L^2}-\frac{\tau}{2}\norm{y^0}^2_{ \mathbb L^2}\right) \right].
\end{equation}
Using the Legendre-Fenchel transform, we obtain that the minimization problem \eqref{Op4} (or equivalently \eqref{Op5})
 can be stated as follow: For every $\tau>0$, find $u^\tau\in\left( L^2(0,T)\right)^{N-1}$, such that
\begin{equation}\label{Eq6}
	J_\tau(u^\tau)=	\inf_{v\in\left(L^2(0,T)\right)^{N-1}} J_\tau(v):=\dis J(v,0)-J(0,0)+\frac{1}{\tau}\norm{I^{1-\gamma}_T \phi^i(0,\cdot;v)}_{ \mathbb L^2}^2 .
\end{equation}
\end{enumerate}
\end{remark}

\begin{definition}
The solution of the minimization problem \eqref{Op4} or equivalently \eqref{Eq6} is called the low-regret control.
\end{definition}

\subsection{Existence and uniqueness of the low regret control}\label{sec51}

We have the following existence and uniqueness result of a low-regret control which is the main result of this section.

\begin{theorem}\label{Existlowregret}
	For every $\tau>0$, there exists a unique $u_\tau\in\left(L^2(0,T)\right)^{N-1}$ solution of the minimization problem \eqref{Eq6}.
\end{theorem}

\begin{proof}
We proceed in several steps.

{\bf Step 1.} First of all, using the fact that $J(v,y^0)\geq 0$ for any $v\in( L^2(0,T))^{N-1}$ and $y^0\in \mathbb L^2$, we can deduce that
	$$J_\tau(v)=\dis J(v,0)-J(0,0)+\frac{1}{\tau}\norm{I^{1-\gamma}_T \phi^i(0,\cdot;v)}_{ \mathbb L^2}^2 \geq -J(0,0).$$ 
	Thus,  $\dis \inf_{v\in\left(L^2(0,T)\right)^{N-1} } J_\tau(v)$ exists. Let $v_n^\tau=(v^{\tau,i}_n)_i$ be a minimizing sequence of $J_\tau$, that is,
	\begin{equation}\label{Conv}
		\lim_{n\to\infty}J_\tau(v_n^\tau)= \inf_{v\in\left(L^2(0,T)\right)^{N-1} } J_\tau(v).
	\end{equation}
	Then, $y_n^\tau:= y(v_n^\tau,0)$ is the very-weak solution  of
\begin{equation}\label{Payvn}
		\left\{
		\begin{array}{lllllllllllllllll}
			\displaystyle	\Dc_t^{\gamma}y_n^{\tau,i}+\mathcal{D}_{b_i^-}^\alpha(\beta^i\mathbb{D}_{a^+}^\alpha y_n^{\tau,i})+q^iy_n^{\tau,i}&=&f^i&\text{in}& Q_i,\,i=1,\dots, N,\\
			\displaystyle	I_{a^+}^{1-\alpha}y_n^{\tau,i}(\cdot,a^+)-I_{a^+}^{1-\alpha}y^{\tau,j}(\cdot,a^+)&=&0&\text{in}&(0,T),~i\neq j=1,\dots, N,\\
			\displaystyle	\sum_{i=1}^N\beta^i(a)\mathbb{D}_{a^+}^\alpha y_n^{\tau,i}(\cdot,a^+)&=&0&\text{in}& (0,T), \\
			\displaystyle	I_{a^+}^{1-\alpha} y^{\tau,1}(\cdot,b_1^-)&=&0&\text{in}&(0,T), \\
			\displaystyle	I_{a^+}^{1-\alpha} y_n^{\tau,i}(\cdot,b_i^-)&=&v^{\tau,i}_n&\text{in}& (0,T), \; i=2,\dots, m\\
			\displaystyle	\beta^i(b_i)\mathbb{D}_{a^+}^\alpha y_n^{\tau,i}(\cdot,b_i^-)&=&v^{\tau,i}_n&\text{in}&(0,T), \; i=m+1,\dots,N,\\
			\displaystyle	y_n^{\tau,i}(0,\cdot) &=&0&\text{in}& (a,b_i),~~i=1,\dots,N.
		\end{array}
		\right.
	\end{equation}
Since $v^{\tau,i}_n\in L^2(0,T),\, ,\, i=2,\cdot, N$ and $f^i\in L^2(Q_i),\, i=1,\cdot, N$, it follows from Theorem \ref{very-weak} that $y_n^\tau\in  L^2((0,T;\mathbb L^2))$. Moreover, there is a constant $C_1>0$ such that
	\begin{equation}\label{Normyn}
		\norm{y_n^\tau}_{L^2((0,T); \mathbb L^2)}\leq C_1 \left(\norm{f}_{L^2((0,T);\mathbb L^2) }+\norm{v_n^\tau}_{\left(L^2(0,T)\right)^{N-1}}\right).
	\end{equation}
	Taking into account \eqref{Conv}, it follows that there is a constant $C>0$ such that
	$$-J(0,0)\leq J_\tau(v_n^\tau)\leq C.$$
	Thus,
	$$0\leq J(v_n^\tau,0)+\frac{1}{\tau}\norm{I^{1-\gamma}_T \phi^i(0,\cdot;v_n)}_{ \mathbb L^2}^2\leq C+J(0,0) \leq C.$$
	Using the form of $J(v_n^\tau,0)$ we obtain that there is a constant $C>0$ such that
	\begin{subequations}
		\begin{alignat}{11}
			\norm{v_n^\tau}_{\left(L^2(0,T)\right)^{N-1}}&\leq C& \label{Normvn}\\
			\norm{I^{1-\gamma}_T \phi(0,\cdot;v_n)}_{ \mathbb L^2}&\leq C&\sqrt{\tau}\label{Normphin}.
		\end{alignat}
	\end{subequations}
	From \eqref{Normyn}, \eqref{Normvn} and \eqref{Payvn}, we obtain that
	\begin{equation}
		\norm{	y_n^\tau}_{L^2((0,T); \mathbb L^2)}\leq C.\label{yn}
	\end{equation}
	Thus,  after subsequences if necessary,  we can deduce that, as $n\to\infty$,
	\begin{subequations}
		\begin{alignat}{5}
			v_n^\tau&\rightharpoonup &u^\tau \,&\text{weakly in}\,\left(L^2(0,T)\right)^{N-1}\label{vnu}\\
			y_n^\tau&\rightharpoonup &y^\tau \,&\text{weakly in}\, L^2((0,T); \mathbb L^2)\label{yny}.
		\end{alignat}
	\end{subequations}
	Since $y_n^\tau$ is the very-weak solution of \eqref{Payvn}, according to Definition \ref{defweak1}, we have that
\begin{equation}\label{sol-n}
		\begin{array}{ll}		&\dis\sum_{i=1}^N\dis\int_{Q_i}y_n^{\tau,i}\left(\D_T^{\gamma}\phi^i+\Dc_{b_i^-}^{\alpha}(\beta^i\Da\phi^i)+q^i\phi^i\right)\,dxdt\\
=&\dis \sum_{i=1}^N\int_{Q_i}f^i\phi^i\,dxdt-\sum_{i=2}^m\int_{0}^Tv^{\tau,i}_n\beta^i(b_i)\Da\phi^i(t,b_i^-)\,dt\\&\dis
+\sum_{i=m+1}^N\int_{0}^Tv_{n}^{\tau,i}(t) I^{1-\alpha}_{a^+}\phi^i(t,b_i^-)\mathrm{d}t,
		\end{array}
	\end{equation}
	for all $\phi=(\phi^i)_i\in\Phi$, where the space $\Phi$ is given in \eqref{defPhi1}.\par
Taking the limit of \eqref{sol-n},  as $n\to\infty$,  while using \eqref{vnu} -\eqref{yny},  we get 
	$$
	\begin{array}{ll}		&\dis\sum_{i=1}^N\dis\int_{Q_i}y^{i,\tau}\left(\D_T^{\gamma}\phi^i+\Dc_{b_i^-}^{\alpha}(\beta^i\Da\phi^i)+q^i\phi^i\right)\,dxdt\\
&=\dis \sum_{i=1}^N\int_{Q_i}f^i\phi^i\,dxdt-\sum_{i=2}^m\int_{0}^Tu^{i,\tau}\beta^i(b_i)\Da\phi^i(t,b_i^-)\,dt\\
&\dis+\sum_{i=m+1}^N\int_{0}^Tu^{i,\tau}(t) I^{1-\alpha}_{a^+}\phi^i(t,b_i^-)\mathrm{d}t.
	\end{array}
	$$
	According to Definition \ref{defweak1},  this implies that $y^\tau=y(u^\tau)$  is the very-weak solution of
	\begin{equation}\label{Paytau}
		\left\{
		\begin{array}{lllllllllllllllll}
			\displaystyle	\Dc_t^{\gamma}y^{i,\tau}+\mathcal{D}_{b_i^-}^\alpha(\beta^i\mathbb{D}_{a^+}^\alpha y^{i,\tau})+q^iy^{i,\tau}&=&f^i&\text{in}& Q_i,\,i=1,\dots, N,\\
			\displaystyle	I_{a^+}^{1-\alpha}y^{i,\tau}(\cdot,a^+)=I_{a^+}^{1-\alpha}y^{j,\tau}(\cdot,a^+)&=&0&\text{in}&(0,T),~i\neq j=1,\dots, N,\\
			\displaystyle	\sum_{i=1}^N\beta^i(a)\mathbb{D}_{a^+}^\alpha y^{i,\tau}(a^+)&=&0&\text{in}& (0,T), \\
			\displaystyle	I_{a^+}^{1-\alpha} y^1(\cdot,b_1^-)&=&0&\text{in}&(0,T), \\
			\displaystyle	I_{a^+}^{1-\alpha} y^{i,\tau}(\cdot,b_i^-)&=&u^{i,\tau}&\text{in}& (0,T), \; i=2,\dots, m\\
			\displaystyle	\beta^i(b_i)\mathbb{D}_{a^+}^\alpha y^{i,\tau}(\cdot,b_i^-)&=&u^{i,\tau}&\text{in}&(0,T), \; i=m+1,\dots,N,\\
			\displaystyle	y^{i,\tau}(0,\cdot;u^\tau) &=&0&\text{in}& (a,b_i),~~i=1,\dots,N.
		\end{array}
		\right.
	\end{equation}
		
	\textbf{Step 2.} Since  $\phi(v_n^\tau)$ is the weak solution of
	\begin{equation}\label{Phivn}
	\left\{\begin{array}{lllllllll}
		\displaystyle	\Dr_T^{\gamma}\phi^i(v_n^\tau)+\mathcal{D}_{b_i^-}^\alpha(\beta^i\mathbb{D}_{a^+}^\alpha \phi^i(v_n^\tau)+q^i\phi^i(v_n^\tau)&=&y^{\tau,i}(v_n^\tau,0)-y^{i}(0,0)&\text{in}& Q_i,\,i=1,\dots, N,\\
		\displaystyle	I_{a^+}^{1-\alpha}\phi^i(a^+;v_n^\tau)-I_{a^+}^{1-\alpha}\phi^j(a^+;v_n^\tau)&=&0&\text{in}&(0,T),~i\neq j=1,\dots, N,\\
		\displaystyle	\sum_{i=1}^N\beta^i(a)\mathbb{D}_{a^+}^\alpha \phi^i(a^+;v_n^\tau)&=&0&\text{in}& (0,T), \\
		\displaystyle	I_{a^+}^{1-\alpha} \phi^i(b_i^-;v_n^\tau)&=&0&\text{in}& (0,T), \; i=1,\dots, m\\
		\displaystyle	\beta^i(b_i)\mathbb{D}_{a^+}^\alpha \phi^i(b_i^-;v_n^\tau)&=&0&\text{in}&(0,T), \; i=m+1,\dots,N,\\
		\displaystyle	I^{1-\gamma}_T \phi^i(T,\cdot;v_n^\tau) &=&0&\text{in}& (a,b_i),~~i=1,\dots,N,
	\end{array}
	\right.
	\end{equation}
	with  $y^{\tau,i}(v_n^\tau,0)-y^{i}(0,0)\in L^2(Q_i),\, i=1,\cdots N$, it follows from Corollary  \ref{coroaux} that
$I_{T}^{1-\gamma} \phi\in C([0,T];\mathbb L^2),\; \phi\in L^2((0,T);\V)~\mbox{ and }~\Dr_T^\gamma \phi\in L^2((0,T);\V^\star).$ Moreover, using \eqref{yn} and \eqref{estaux},we have that there is a constant $C>0$ independent of $n$ such that
\begin{subequations}
\begin{alignat}{5}
			\norm{\phi(v_n^\tau)}_{L^2((0,T);\V)}&\leq& C\label{mwj0},\\
			\norm{I^{1-\gamma}_{T}(\phi(v_n^\tau))}_{C([0,T];\mathbb{L}^2)}&\leq& C\label{mwj1},\\
			\norm{\Dr_T^{\gamma}\phi^i(v_n^\tau)}_{L^2((0,T);\V^\star)}&\leq& C\label{mwj2}.
		\end{alignat}
	\end{subequations}
From \eqref{mwj0}, \eqref{mwj2} and Lemma \ref{lemK}, we can deduce that
\begin{equation}\label{mwj1K}
 \norm{\phi(v_n^\tau)}_{\mathbb{K}}\leq C,
\end{equation}
for some $C>0$ independent of $n$, where $\mathbb{K}$ is the Hilbert space given by \eqref{defK}.  Hence, there exists $\phi^\tau\in \mathbb{K}$ such that as $n\to\infty$,
\begin{subequations}
\begin{alignat}{11}
  \phi(v_n^\tau)\rightharpoonup \phi^\tau\,\,\text{weakly in }\,\, L^2((0,T);\V)\label{convphin}\\
  D_T^{\gamma}\phi(v_n^\tau)\rightharpoonup D_T^{\gamma}\phi^\tau\,\,\text{weakly in } L^2((0,T);\V^\star)\label{convDphin}.
\end{alignat}
\end{subequations}
Moreover, using Remark \ref{rem28} and the fact that $\V\hookrightarrow \mathbb{L}^2\hookrightarrow \V^\star$,  it follows from \eqref{contWTA} that
 \begin{equation}\label{IT0}I^{1-\gamma}_{T}\phi^\tau\in C([0,T],\mathbb{L}^2),
 \end{equation}
 and  from the last equation in \eqref{Phivn} we have that
 \begin{equation}\label{IT}
 I^{1-\gamma}_{T}\phi^\tau(T,\cdot)=0.
\end{equation}

On the other hand, $\phi(v_n^\tau)\in L^2((0,T);\V)$ being  the weak solution of \eqref{Phivn}, we have from Definition \ref{definition47} that for every $\varphi\in \V$ and a.e. $t\in (0,T)$,
$$\dis \langle \Dr_T^\gamma \phi(v_n^\tau)(t,\cdot), \varphi\rangle_{\V^\star, \V}+\mathbb F(\phi(v_n^\tau)(t,\cdot),\varphi)=\langle y^{\tau}(v_n^\tau,0)(t,\cdot)-y(0,0)(t,\cdot),\varphi\rangle_{\V^\star, \V}. $$
Passing to the limit as $n\to \infty$ in this latter identity while using \eqref{yny} and \eqref{convphin},we obtain that
for every $\varphi\in \V$ and a.e.  $t\in (0,T)$, 
\begin{equation}\label{IT1}
\dis \langle \Dr_T^\gamma \phi^\tau(t,\cdot), \varphi\rangle_{\V^\star, \V}+\mathbb F(\phi^\tau(t,\cdot),\varphi)=\langle y^{\tau}(t,\cdot)-y(0,0)(t,\cdot),\varphi\rangle_{\V^\star, \V}.
\end{equation}
From \eqref{IT0}, \eqref{convphin}, \eqref{convDphin}, \eqref{IT} and \eqref{IT1}, we have that $\phi^\tau=\phi(u^\tau)$ is a weak solution of \eqref{PaPhi} with $v=u^\tau$. Thus $\phi^\tau=\phi(u^\tau)$ satisfies
\begin{equation}\label{phitau}\dis
		\left\{\begin{array}{lllllllll}
			\mathbb{D}_T^\gamma \phi^{i,\tau}+\dis  \Dc_{b_i^-}^\alpha\,(\beta^i\, \Dr_{a^+}^\alpha \phi^{i,\tau})+q\,\phi^{i,\tau}&=&y^{i,\tau}(u^\tau,0)-y^{i}(0,0)&\hbox{ in }& Q_i,i=1,\dots, N\\
			\dis (I_{a^+}^{1-\alpha} \phi^{i,\tau})(\cdot,a^+)- (I_{a^+}^{1-\alpha} \phi^{j,\tau})(\cdot,a^+)&=&0&\hbox{ in }& (0,T),i,j=1,\dots, N\\\dis \sum_{i=1}^N\beta^i(a)\D_a^\alpha\phi^{i,\tau}(\cdot,a^+)&=&0&\hbox{in}& (0,T),\\	I_{a^+}^{1-\alpha} \phi^{i,\tau}(\cdot,b_i^-)&=&0&\text{in}& (0,T), \; i=1,\dots, m\\
			\dis (\beta^i\Dr_{a^+}^{\alpha}\phi^{i,\tau})(\cdot,b_i^-)&=&0&\hbox{ in }& (0,T),i=m+1,\cdots,N\\
			I^{1-\gamma}_T\phi^{i,\tau}(T,\cdot;u^\tau)&=&0&\hbox{ in } &(a,b_i), i=1,\cdots,N.
		\end{array}
		\right.
	\end{equation}
Using again \eqref{IT0} and \eqref{Normphin} we get that as $n\to\infty$,
\begin{equation}\label{convphi0n}
 I^{1-\gamma}_T \phi^i(0,\cdot;v_n)\rightharpoonup I^{1-\gamma}_T\phi^{i,\tau}(0,\cdot;u^\tau) \hbox{ weakly in } \mathbb{L}^2.
\end{equation}

\textbf{Step 3:} Since the functional $v\mapsto J_\tau(v)$ is continuous, it follows that it is lower semi-continuous. Thus, using \eqref{vnu}, \eqref{yny} and \eqref{convphi0n}, we deduce that
	$$\dis J_\tau(u_\tau)\leq\liminf_{n\to+\infty}J_\tau (v_n),$$
	which combining with \eqref{Conv} implies that 
	\begin{equation*}
		J_\tau(u^\tau)=\inf_{v\in\left(L^2(0,T)\right)^{N-1}}J_\tau(v).
	\end{equation*} 
	The uniqueness of the low-regret control $u^\tau$ follows from the strict convexity of the functional $J_\tau$. The proof is finished.
\end{proof}

\subsection{Characterization of the low-regret control }
In this section we give a complete characterization (optimality systems and optimality conditions) of the low-regret control. The following is the main result of this section.

\begin{theorem}
Let $u^\tau\in \left(L^2(0,T)\right)^{N-1}$ be the low-regret control constructed in Theorem \ref{Existlowregret}.  Then,  there exist $h^\tau\in L^2((0,T); \V)$, $p^\tau\in L^2((0,T); \V)$ such that $(u^\tau, y^{\tau},\phi^{\tau}, h^\tau, p^\tau)$ satisfies the following.  We have that $y^\tau\in L^2((0,T); \mathbb L^2)$ is the very-weak solution of \eqref{Paytau}, $\phi^{\tau}\in L^2((0,T); \V)$ is the weak solution of \eqref{phitau},
$h^\tau$ is the weak solution of	
\begin{equation}
			\left\{
				\begin{array}{lllllllllllllllll}
					\displaystyle	\Dc_t^{\gamma}h^{i,\tau}+\mathcal{D}_{b_i^-}^\alpha(\beta^i\mathbb{D}_{a^+}^\alpha h^{i,\tau})+q^ih^{i,\tau}&=&0&\text{in}& Q_i,\,i=1,\dots, N,\\
					\displaystyle	I_{a^+}^{1-\alpha}h^{i,\tau}(\cdot,a^+)-I_{a^+}^{1-\alpha}h^{j,\tau}(\cdot,a^+)&=&0&\text{in}&(0,T),~i\neq j=1,\dots, N,\\
					\displaystyle	\sum_{i=1}^N\beta^i(a)\mathbb{D}_{a^+}^\alpha h^{i,\tau}(\cdot,a^+)&=&0&\text{in}& (0,T), \\
					\displaystyle	I_{a^+}^{1-\alpha} h^{i,\tau}(\cdot,b_i^-)&=&0&\text{in}& (0,T), \; i=1,\dots, m\\
					\displaystyle	\beta^i(b_i)\mathbb{D}_{a^+}^\alpha h^{i,\tau}(\cdot,b_i^-)&=&0&\text{in}&(0,T), \; i=m+1,\dots,N,\\
					\displaystyle	h^{i,\tau}(0,\cdot) &=&\frac{1}{\sqrt{\tau}}I^{1-\gamma}_T\phi^{i,\tau}(0,\cdot;u^\tau)&\text{in}& (a,b_i),~~i=1,\dots,N,
				\end{array}
				\right.\label{paqtau}
		\end{equation}
and $p^\tau$ is the weak solution of
\begin{equation}
		\left\{
		\begin{array}{lllllllllllllllll}
			\displaystyle			\mathbb{D}_T^\gamma p^{i,\tau}+\mathcal{D}_{b_i^-}^\alpha(\beta^i\mathbb{D}_{a^+}^\alpha p^{i,\tau})+q^ip^{i,\tau}&=&y^{i,\tau}-y_d^i+\frac{1}{\sqrt{\tau}}h^{i, \tau}&\text{in}& Q_i,\,i=1,\dots, N,\\
			\displaystyle	I_{a^+}^{1-\alpha}p^{i,\tau}(\cdot,a^+)-I_{a^+}^{1-\alpha}p^{j,\tau}(\cdot,a^+)&=&0&\text{in}&(0,T),~i\neq j=1,\dots, N,\\
			\displaystyle	\sum_{i=1}^N\beta^i(a)\mathbb{D}_{a^+}^\alpha p^{i,\tau}(\cdot,a^+)&=&0&\text{in}& (0,T), \\
			\displaystyle	I_{a^+}^{1-\alpha} p^{i,\tau}(\cdot,b_i^-)&=&0&\text{in}& (0,T), \; i=1,\dots, m\\
			\displaystyle	\beta^i(b_i)\mathbb{D}_{a^+}^\alpha p^{i,\tau}(\cdot,b_i^-)&=&0&\text{in}&(0,T), \; i=m+1,\dots,N,\\
			\displaystyle	I^{1-\gamma}_Tp^{i,\tau}(T,\cdot) &=&0&\text{in}& (a,b_i),~~i=1,\dots,N.
		\end{array}
	\right.\label{paptau}
	\end{equation}
In addition, we have that for all $v\in\left (L^2(0,T)\right)^{N-1}$,
\begin{equation}\label{inter7}
\begin{array}{rll}
\dis \sum_{i=2}^N\int_0^T \Big(\zeta u^{i,\tau}-\beta^i(b_i)\mathbb{D}_{a^+}^\alpha p^{i,\tau}(\cdot,b_i^-)+I^{1-\alpha}_{a^+}p^{i,\tau}(\cdot,b_i^-)\Big)(v^i-u^{i,\tau})dt=0.
\end{array}
		\end{equation}
\end{theorem}

\begin{proof}
We have already shown that $y^\tau$ satisfies  \eqref{Paytau}.
	Writing the Euler-Lagrange first order optimality conditions  for the low-regret control $u^\tau$, we get for all $v\in\left(L^2(0,T)\right)^{N-1}$,
	\begin{equation}\label{euler1}
		\lim_{\mu\downarrow 0}\frac{J_\tau(u^\tau+\mu(v-u^\tau))-J_\tau(u^\tau)}{^\mu}=0.
	\end{equation}
After some computations, \eqref{euler1} gives for all $v\in \left(L^2(0,T)\right)^{N-1}$,
\begin{align}\label{inter8}
		\sum_{i=1}^N\int_{Q_i} z^i(v-u^\tau)(y^{i}(u^\tau,0)-y_d^i)\,dxdt &+\sum_{i=1}^N\int_a^{b_i} \frac{1}{\sqrt{\tau}}I^{1-\gamma}_T\phi_z^i(0,\cdot;v-u^\tau)\frac{1}{\sqrt{\tau}}I^{1-\gamma}_T\phi^i(0,\cdot;u^\tau)\,dx\notag\\
		&+\zeta\sum_{i=2}^N\int_0^T u^{i,\tau} (v^i-u^{i ,\tau})\,dt=0, 
	\end{align}
where  $\phi^\tau$ is  the weak solution of  \eqref{phitau},  $\dis z^i(v-u^\tau):=\frac{y^{i}(u^\tau+\mu(v-u^\tau),0)-y^{i}(u^\tau,0)}{\mu}\in L^2((0,T);\mathbb{L}^2)$ is the very-weak solution of
	\begin{equation}\label{Pazv}
		\left\{\begin{array}{lllllllll}\Dc_t^{\gamma}z^i(v-u^\tau)+\mathcal{D}_{b_i^-}^\alpha(\beta^i\mathbb{D}_{a^+}^\alpha z^i(v-u^\tau))+q^iz^i(v-u^\tau)&=&0&\text{in}& Q_i,\,i=1,\dots, N,\\
		\displaystyle	I_{a^+}^{1-\alpha}z^i(a^+;v-u^\tau)-I_{a^+}^{1-\alpha}z^{j}(a^+;v-u^\tau)&=&0&\text{in}&(0,T),~i\neq j=1,\dots, N,\\
		\displaystyle	\sum_{i=1}^N\beta^i(a)\mathbb{D}_{a^+}^\alpha z^i(a^+;v-u^\tau)&=&0&\text{in}& (0,T), \\
		\displaystyle	I_{a^+}^{1-\alpha} z^{1}(b_1^-;v-u^\tau)&=&0&\text{in}&(0,T), \\
		\displaystyle	I_{a^+}^{1-\alpha} z^i(b_i^-;v-u^\tau)&=&v^{i}-u^{i,\tau}&\text{in}& (0,T), \; i=2,\dots, m\\
		\displaystyle	\beta^i(b_i)\mathbb{D}_{a^+}^\alpha z^i(b_i^-;v-u^\tau)&=&v^{i}-u^{i,\tau}&\text{in}&(0,T), \; i=m+1,\dots,N,\\
		\displaystyle	z^i(0;v-u^\tau) &=&0&\text{in}& (a,b_i),~~i=1,\dots,N
		\end{array}
		\right.
	\end{equation}
	and  $\dis \phi_z(v-u^\tau):=\frac{\phi(u^\tau+\mu(v-u^\tau))-\phi(u^\tau)}{\mu}\in L^2((0,T);\V)$ is the weak solution of
	\begin{equation}\label{phiztau}
		\left\{\begin{array}{lllllllllllllllll}
			\displaystyle		\mathbb{D}_T^\gamma \phi_z^i(v-u^\tau)+\mathcal{D}_{b_i^-}^\alpha(\beta^i\mathbb{D}_{a^+}^\alpha \phi_z^i(v-u^\tau))+q^i\phi_z^i(v-u^\tau)&=&z(v-u^\tau)&\text{in}& Q_i,\,i=1,\dots, N,\\
			\displaystyle	I_{a^+}^{1-\alpha}\phi_z^i(a^+;v-u^\tau)-I_{a^+}^{1-\alpha}\phi_z^j(a^+;v-u^\tau)&=&0&\text{in}&(0,T),~i\neq j=1,\dots, N,\\
			\displaystyle	\sum_{i=1}^N\beta^i(a)\mathbb{D}_{a^+}^\alpha \phi_z^i(a^+;v-u^\tau)&=&0&\text{in}& (0,T), \\
			\displaystyle	I_{a^+}^{1-\alpha} \phi_z^i(b_i^-;v-u^\tau)&=&0&\text{in}& (0,T), \; i=1,\dots, m\\
			\displaystyle	\beta^i(b_i)\mathbb{D}_{a^+}^\alpha \phi_z^i(b_i^-;v-u^\tau)&=&0&\text{in}&(0,T), \; i=m+1,\dots,N,\\
			\displaystyle	I^{1-\gamma}_T\phi_z^i(T;v-u^\tau) &=&0&\text{in}& (a,b_i),~~i=1,\dots,N.
		\end{array}
		\right.
	\end{equation}
In view of \eqref{convphi0n}, $ I^{1-\gamma}_T\phi^{\tau}(0,\cdot,u^\tau) \in  \mathbb{L}^2.$ Thus, from Theorem \ref{Theo1aux1} we have that $h^{\tau}\in L^2((0,T);\V)$  is the unique weak solution of \eqref{paqtau}. Moreover, there is a constant $C>0$ such that
\begin{equation}\label{sam3}
	\norm{h^{\tau}}_{L^2(0,T,\V)} +\norm{I^{1-\gamma}_{t}h^{\tau}}_{C([0,T];\mathbb L^2)}+	\norm{\Dr_t^{\gamma}h^{\tau}}_{L^2(0,T;\V^\star)}\leq C
\frac{1}{\sqrt{\tau}}\norm{I^{1-\gamma}_T\phi^{i,\tau}(0,\cdot;u^\tau)}_{\mathbb{L}^2}.
\end{equation}
Thus,  if we multiply the first equation in  \eqref{phiztau} by $\dis \frac{1}{\sqrt{\tau}}h^{i, \tau}$ and  integrate by parts over $Q_i$ we get
\begin{equation}\label{inter9}
		\sum_{i=1}^N\int_{Q_i}z^i(v-u^\tau) \frac{1}{\sqrt{\tau}}h^{i, \tau}\,dxdt=\sum_{i=1}^N\int_{a}^{b_i} \frac{1}{\sqrt{\tau}}I^{1-\gamma}_T\phi_z^i(0,\cdot;v-u^\tau)\frac{1}{\sqrt{\tau}}I^{1-\gamma}_T\phi^{i,\tau}(0,\cdot;u^\tau)\,dx.
	\end{equation}
	Combining \eqref{inter8}-\eqref{inter9} gives for all $v\in\left(L^2(0,T)\right)^{N-1}$,
\begin{equation}\label{inter10}
		\sum_{i=1}^N\int_Q z^i(v-u^\tau)\left[y^{i}(u^\tau,0)-y_d^i+ \frac{1}{\sqrt{\tau}}h^{i, \tau}\right]\,dxdt+\sum_{i=2}^N\int_0^T \zeta u^{i,\tau} (v^i-u^{i,\tau})\,dt=0.
	\end{equation}
	Since $y(u^\tau,0)-y_d+ \frac{1}{\sqrt{\tau}}h^{\tau}\in L^2((0,T); \mathbb L^2)$, then according to Corollary  \ref{coroaux}, the weak solution $p^\tau$ of \eqref{paptau} satisfies $p^\tau\in L^2((0,T);\V)$.  Moreover, there is a constant $C>0$ such that \begin{equation}\label{sam4}
	\norm{p^{\tau}}_{L^2(0,T,\V)} +\norm{I^{1-\gamma}_{T}p^{\tau}}_{C([0,T];\mathbb L^2)}+	\norm{\Dr_T^{\gamma}p^{\tau}}_{L^2(0,T;\V^\star)}\leq C \norm{y^{\tau}-y_d+\frac{1}{\sqrt{\tau}}h^{\tau}}_{L^2(0,T;\mathbb{L}^2)}.
\end{equation}
It  follows from  \eqref{traceI} and \eqref{TraceD}  that $I^{1-\alpha}_{a^+}p^{i,\tau}(\cdot,b_i^-)\in L^2(0,T)$ and $\mathbb{D}_{a^+}^\alpha p^{i,\tau}(\cdot,b_i^-)\in L^2(0,T),$ $i=1,\cdots, N$. Thus,  multiplying \eqref{Pazv} by $p^\tau$ and integrating by parts, we obtain that
	\begin{align}\label{inter11}
&\dis	\sum_{i=1}^N	\int_{Q_i} z^i(v-u^\tau)\left[y^{i}(u^\tau,0)-y_d^i+ \frac{1}{\sqrt{\tau}}h^{i, \tau}\right]\,dxdt\notag\\
=&
\dis -\sum_{i=2}^m\int_0^T \beta^i(b_i)\mathbb{D}_{a^+}^\alpha p^{i,\tau}(\cdot,b_i^-)(v^i-u^{i,\tau})dt
+\dis
\sum_{i=m+1}^N\int_0^T(v^i-u^{i,\tau})I^{1-\alpha}_{a^+}p^{i,\tau}(t,b_i^-)dt.
	\end{align}
Combining \eqref{inter10}-\eqref{inter11}, we get for all $v\in\left(L^2(0,T)\right)^{N-1}$,
$$
	\dis \sum_{i=2}^N\int_0^T \left(\zeta u^{i,\tau}-\beta^i(b_i)\mathbb{D}_{a^+}^\alpha p^{i,\tau}(\cdot,b_i^-)+I^{1-\alpha}_{a^+}p^{i,\tau}(\cdot,b_i^-)\right)(v^i-u^{i,\tau})dt=0, 
$$
from which we can deduce \eqref{inter7}.The proof is finished.
\end{proof}

\subsection{Existence and characterization of the no-regret control}
In this section we prove the existence, uniqueness and we characterize the no-regret control. We start with the  existence  and uniqueness result.

\begin{proposition}\label{prop-58}
The low-regret control $u^\tau$ constructed in Theorem \ref{Existlowregret} converges, as $\tau\downarrow 0$, to a no-regret control $u$ which is the solution of the minimization problem \eqref{inte4}.
\end{proposition}

\begin{proof} 
Since $u^\tau$ is the low-regret control,  we can deduce that $J_\tau(u^\tau)\leq J_\tau(v)$ for any $v\in \left(L^2(0,T)\right)^{N-1}.$ Consequently,
$$J_\tau(u^\tau)\leq  J_\tau(0).$$
Therefore, using the definition of $J_\tau$ given by \eqref{Eq6}, we have that
\begin{equation}\label{sam2}J(u^\tau ,0)-J(0,0)+\frac{1}{\tau}\norm{I^{1-\gamma}_T \phi^i(0,\cdot,u^\tau)}^2_{ \mathbb L^2}\leq \dis\frac{1}{\tau}\norm{I^{1-\gamma}_T \phi^i(0,\cdot,0)}^2_{ \mathbb L^2}.\end{equation}
Observing that $\phi:=\phi(v)$ is the weak solution of \eqref{PaPhi} with $v=0$, we obtain from \eqref{sam1} that
$\phi^i(0)=0$ on $Q_i,\,1\leq i\leq N$ and \eqref{sam2} becomes
$$J(u^\tau ,0)-J(0,0)+\frac{1}{\tau}\norm{I^{1-\gamma}_T \phi^i(0,\cdot,u^\tau)}^2_{ \mathbb L^2}\leq 0,$$
 which in view of the  expression of $J$ given in \eqref{mp2} implies that
%
\begin{align}\label{8Inter}
		\norm{y(u^\tau,0)-y_d}^2_{L^2((0,T); \mathbb L^2)}&+\zeta\norm{u^\tau}^2_{\left(L^2(0,T)\right)^{N-1}}+\frac{1}{\tau}\norm{I^{1-\gamma}_T \phi^i(0,\cdot,u^\tau)}^2_{ \mathbb L^2}\notag\\
&\leq \norm{y(0,0)-y_d}^2_{L^2((0,T); \mathbb L^2)}.
	\end{align}
Consequently,
	\begin{subequations}
		\begin{alignat}{5}
			&\norm{y(u^\tau,0)}_{L^2((0,T); \mathbb L^2)}&\leq & \norm{y(0,0)}_{L^2((0,T); \mathbb L^2)}+\norm{y_d}_{L^2((0,T); \mathbb L^2)},\label{1inter}\\
&\norm{u^\tau}_{\left(L^2(0,T)\right)^{N-1}}&\leq &\norm{y(0,0)-y_d}_{L^2((0,T); \mathbb L^2)},\label{2inter}\\
&\norm{I^{1-\gamma}_T \phi^i(0,\cdot;u^\tau)}_{ \mathbb L^2}&\leq& \sqrt{\tau}\norm{y(0,0)-y_d}_{L^2((0,T); \mathbb L^2)}\label{3inter}.
		\end{alignat}
	\end{subequations}
Since $\phi^\tau:=\phi(u^\tau)$ is the weak solution of \eqref{PaPhi} with $v=u^\tau$, using \eqref{1inter}, we get from \eqref{sam1} that there is a constant $C>0$ such that
\begin{subequations}\label{sam1bis}
		\begin{alignat}{5}
			&\norm{	\Dr_T^{\gamma}\phi^\tau }_{L^2((0,T);\V^\star)}&\leq & C\left(\norm{y(0,0)}_{L^2((0,T);\mathbb L^2)}+\norm{y_d}_{L^2((0,T); \mathbb L^2)}\right),\label{sam1a}\\
&\norm{\phi^\tau }_{L^2(0,T,\V)}&\leq &C\left(\norm{y(0,0)}_{L^2((0,T);\mathbb L^2)}+\norm{y_d}_{L^2((0,T); \mathbb L^2)}\right),\label{sam1b}\\
&\norm{I^{1-\gamma}_T\phi}_{C([0,T];\mathbb L^2)}&\leq& C\left(\norm{y(0,0)}_{L^2((0,T);\mathbb L^2)}+\norm{y_d}_{L^2((0,T) \mathbb L^2)}\right).\label{sam1c}
		\end{alignat}
	\end{subequations}
From \eqref{convphi0n}, we have that $I^{1-\gamma}_T\phi^{\tau}(0,\cdot,u^\tau)\in \mathbb{L}^2$. It then follows from Theorem \ref{Theo1aux1} that there exists a unique  $h^{\tau}\in L^2(0,T;\V)$ weak solution of \eqref{paqtau}. Moreover, there is a constant $C>0$ such that $$\begin{array}{llll}
			\norm{h^{\tau}}_{L^2(0,T;\V)}&\leq & \dis \frac{C}{\sqrt{\tau}}\norm{I^{1-\gamma}_T\phi^{i,\tau}(0,\cdot;u^\tau)}_{\mathbb L^2},\\
\norm{I^{1-\gamma}_{t}h^{\tau}}_{C([0,T];\mathbb L^2)}&\leq &\dis \frac{C}{\sqrt{\tau}}\norm{I^{1-\gamma}_T\phi^{i,\tau}(0,\cdot;u^\tau)}_{\mathbb L^2}\\
\norm{\Dr_t^{\gamma}h^{\tau}}_{L^2(0,T;\V^\star)}&\leq&\dis \frac{C}{\sqrt{\tau}}\norm{I^{1-\gamma}_T\phi^{i,\tau}(0,\cdot;u^\tau)}_{\mathbb L^2},
	\end{array}
$$
which according to \eqref{3inter} implies that
\begin{subequations}\label{sam5bis}
		\begin{alignat}{5}
			&\norm{h^{\tau}}_{L^2(0,T;\V)}&\leq & C\left(\norm{y(0,0)}_{L^2((0,T);\mathbb L^2)}+\norm{y_d}_{L^2((0,T); \mathbb L^2)}\right),\label{sam5a}\\
&\norm{I^{1-\gamma}_{t}h^{\tau}}_{C([0,T];\mathbb L^2)}&\leq &C\left(\norm{y(0,0)}_{L^2((0,T);\mathbb L^2)}+\norm{y_d}_{L^2((0,T); \mathbb L^2)}\right),\label{sam5b}\\
&\norm{\Dr_t^{\gamma}h^{\tau}}_{L^2(0,T;\V^\star)}&\leq& C\left(\norm{y(0,0)}_{L^2((0,T);\mathbb L^2)}+\norm{y_d}_{L^2((0,T); \mathbb L^2)}\right),\label{sam5c}
		\end{alignat}
	\end{subequations}
for some constant $C>0$ independent of $\tau$.

From \eqref{2inter}, \eqref{1inter} and  \eqref{sam1b}, we can deduce that there exist $u\in \left(L^2(0,T)\right)^{N-1},$ $y\in L^2((0,T); \mathbb L^2)$  and $\phi\in L^2(0,T;\V)$  such that as $\tau\to 0$,
\begin{subequations}\label{12interbis}
		\begin{alignat}{7}
&{u^\tau}&\rightharpoonup& u&\text{weakly in}\,\,& {\left(L^2(0,T)\right)^{N-1}},\label{12inter}\\
			&{y(u^\tau,0)}&\rightharpoonup& y\,\,&\text{weakly in}\,& L^2((0,T); \mathbb L^2) \label{12intery} ,\\
&{\phi(u^\tau)}&\rightharpoonup& \phi\,\,&\text{weakly in}\,& L^2((0,T); \V) \label{12interphi}.
		\end{alignat}
	\end{subequations}
Therefore, proceeding as for $y_n^\tau$ and $\phi(v_n^\tau)$ in Section \ref{sec51}  while using \eqref{sam1bis} and \eqref{12interbis}, we obtain that
$y=y(u,0)\in L^2((0,T);\mathbb{L}^2)$ and  $\phi=\phi(u)\in L^2((0,T);\V)$  are respectively solutions of
	\begin{equation}\label{pay}
	\left\{
	\begin{array}{lllllllllllllllll}
		\displaystyle	\Dc_t^{\gamma}y^i+\mathcal{D}_{b_i^-}^\alpha(\beta^i\mathbb{D}_{a^+}^\alpha y^{i})+q^iy^{i}&=&f^i&\text{in}& Q_i,\,i=1,\dots, N,\\
		\displaystyle	I_{a^+}^{1-\alpha}y^{i}(\cdot,a^+)-I_{a^+}^{1-\alpha}y^{j}(\cdot,a^+)&=&0&\text{in}&(0,T),~i\neq j=1,\dots, N,\\
		\displaystyle	\sum_{i=1}^N\beta^i(a)\mathbb{D}_{a^+}^\alpha y^{i}(\cdot,a^+)&=&0&\text{in}& (0,T), \\
		\displaystyle	I_{a^+}^{1-\alpha} y^1(\cdot,b_1^-)&=&0&\text{in}&(0,T), \\
		\displaystyle	I_{a^+}^{1-\alpha} y^{i}(\cdot,b_i^-)&=&u^i&\text{in}& (0,T), \; i=2,\dots, m\\
		\displaystyle	\beta^i(b_i)\mathbb{D}_{a^+}^\alpha y^{i}(\cdot,b_i^-)&=&u^i&\text{in}&(0,T), \; i=m+1,\dots,N,\\
		\displaystyle	y^{i}(0,\cdot) &=&0&\text{in}& (a,b_i),~~i=1,\dots,N,
	\end{array}
	\right.
	\end{equation}	
and
	\begin{equation}\label{phiu}
			\left\{\begin{array}{lllllllllllllllll}
			\displaystyle		\mathbb{D}_T^\gamma \phi^i(u)+\mathcal{D}_{b_i^-}^\alpha(\beta^i\mathbb{D}_{a^+}^\alpha \phi^i(u))+q^i\phi^i(u)&=&y^i(u,0)-y^i(0,0)&\text{in}& Q_i,\,i=1,\dots, N,\\
			\displaystyle	I_{a^+}^{1-\alpha}\phi^i(u)(\cdot,a^+)-I_{a^+}^{1-\alpha}\phi^j(u)(\cdot,a^+)&=&0&\text{in}&(0,T),~i\neq j=1,\dots, N,\\
			\displaystyle	\sum_{i=1}^N\beta^i(a)\mathbb{D}_{a^+}^\alpha \phi^i(u)(\cdot,a^+)&=&0&\text{in}& (0,T), \\
			\displaystyle	I_{a^+}^{1-\alpha} \phi^i(u)(\cdot,b_i^-)&=&0&\text{in}& (0,T), \; i=1,\dots, m\\
			\displaystyle	\beta^i(b_i)\mathbb{D}_{a^+}^\alpha \phi^i(u)(\cdot,b_i^-)&=&0&\text{in}&(0,T), \; i=m+1,\dots,N,\\
			\displaystyle	I^{1-\gamma}_T\phi^i(u)(T,\cdot) &=&0&\text{in}& (a,b_i),~~i=1,\dots,N.
		\end{array}
		\right.
	\end{equation}
Since $y^i(u,0)-y^i(0,0)\in L^2(Q_i),\,\, i=1,\cdots N$, it follows from Corollary  \ref{coroaux} that
$\phi(u)\in L^2((0,T);\V)$ and $I^{1-\gamma}_{T}(\phi(u))\in C([0,T];L^2(\Omega))$. Thus, $I^{1-\gamma}_{T}(\phi(0,\cdot;u))\in L^2(\Omega)$ and in view of \eqref{3inter}, we can deduce that
	\begin{equation}\label{4inter}
		I^{1-\gamma}_{T}(\phi(0,\cdot;u^\tau)\to 	I^{1-\gamma}_{T}(\phi(0,\cdot;u))=0\,\,\text{strongly in }\, \mathbb L^2\;\mbox{ as } \tau\to 0.
	\end{equation}
	Consequently, $( y^0,I^{1-\gamma}_{T}(\phi(0,\cdot,u)))_{\mathbb L^2}=0$, which means that $u$ belongs to the set  $ \mathcal{U}$ defined in \eqref{defU}. In other words $u\in \mathcal{U}$ is the unique solution of the no-regret control problem \eqref{inte4}. The proof is finished.
	\end{proof}

Next, we characterize the no-regret control constructed in Proposition \ref{prop-58}.
	
	\begin{theorem}
Let $\mathcal{U}$ be defined as in \eqref{defU} and $u\in \mathcal{U}$ be the no-regret control associated to the state $y=y(u,0)$. Then,  there exist $h\in L^2((0,T); \V)$, $p\in L^2((0,T);\V)$, $\zeta_1,\zeta_2\in L^2((0,T); \mathbb L^2)$  such that $y\in L^2((0,T); \mathbb L^2)$ is the very-weak solution of \eqref{pay}, $h$ is the weak solution of
		\begin{equation}		
			\left\{
		\begin{array}{lllllllllllllllll}
		\displaystyle	\Dc_t^{\gamma}h^i+\mathcal{D}_{b_i^-}^\alpha(\beta^i\mathbb{D}_{a^+}^\alpha h^i)+q^ih^i&=&0&\text{in}& Q_i,\,i=1,\dots, N,\\
		\displaystyle	I_{a^+}^{1-\alpha}h^i(\cdot,a^+)-I_{a^+}^{1-\alpha}h^{j}(\cdot,a^+)&=&0&\text{in}&(0,T),~i\neq j=1,\dots, N,\\
		\displaystyle	\sum_{i=1}^N\beta^i(a)\mathbb{D}_{a^+}^\alpha h^i(\cdot,a^+)&=&0&\text{in}& (0,T), \\
		\displaystyle	I_{a^+}^{1-\alpha} h^i(\cdot,b_i^-)&=&0&\text{in}& (0,T), \; i=1,\dots, m\\
		\displaystyle	\beta^i(b_i)\mathbb{D}_{a^+}^\alpha h^i(\cdot,b_i^-)&=&0&\text{in}&(0,T), \; i=m+1,\dots,N,\\
	\displaystyle	h^i(0,\cdot) &=&\zeta^i_1&\text{in}& (a,b_i),~~i=1,\dots,N,
		\end{array}
		\right.\label{paq}
	\end{equation}
	and $p$ is the weak solution of
	\begin{equation}
		\left\{
		\begin{array}{lllllllllllllllll}
	\displaystyle	\mathbb{D}_T^\gamma p^{i,\tau}+\mathcal{D}_{b_i^-}^\alpha(\beta^i\mathbb{D}_{a^+}^\alpha p^i)+q^ip^i&=&y^{i}-y_d^i+\zeta^i_2&\text{in}& Q_i,\,i=1,\dots, N,\\
	\displaystyle	I_{a^+}^{1-\alpha}p^i(\cdot,a^+)-I_{a^+}^{1-\alpha}p^{j}(\cdot,a^+)&=&0&\text{in}&(0,T),~i\neq j=1,\dots, N,\\
	\displaystyle	\sum_{i=1}^N\beta^i(a)\mathbb{D}_{a^+}^\alpha p^i(\cdot,a^+)&=&0&\text{in}& (0,T), \\
	\displaystyle	I_{a^+}^{1-\alpha} p^i(\cdot,b_i^-)&=&0&\text{in}& (0,T), \; i=1,\dots, m\\
	\displaystyle	\beta^i(b_i)\mathbb{D}_{a^+}^\alpha p^i(\cdot,b_i^-)&=&0&\text{in}&(0,T), \; i=m+1,\dots,N,\\
	\displaystyle	I^{1-\gamma}_Tp^i(T,\cdot) &=&0&\text{in}& (a,b_i),~~i=1,\dots,N.
		\end{array}
		\right.\label{pap}
	\end{equation}
In addition, we have that for all $v\in\left(L^2(0,T)\right)^{N-1}$,
	\begin{equation}
\sum_{i=2}^N\int_0^T(v^i-u^{i})\left[ u^{i} +I^{1-\alpha}_{a^+}p^{i}(t,b_i^-)-\beta^i(b_i)\mathbb{D}_{a^+}^\alpha p^i(\cdot,b_i^-)\right]dt=0.\label{7inter}
		\end{equation}
	\end{theorem}
	
	\begin{proof}
		We have already shown that $y$ is the very-weak solution of \eqref{pay}. From \eqref{3inter}, we obtain that there is a constant $C>0$ such that
		\begin{equation*}
			\norm{\frac{1}{\sqrt{\tau}}\phi(u^\tau)}_{L^2((0,T); \mathbb L^2)}\leq C.
		\end{equation*}
		Therefore, it follows from  \eqref{estim_y1} that
		$$\norm{h^{ \tau}}_{L^2((0,T);\V)}\leq C_1	\norm{\frac{1}{\sqrt{\tau}}\phi(u^\tau)}_{L^2((0,T); \mathbb L^2)}\leq C.$$
		Using \eqref{inter10},   we get that
		\begin{equation*}\begin{array}{llll}
		&\dis \left(z(v-u^\tau), \frac{1}{\sqrt{\tau}}h^{ \tau}\right)_{L^2((0,T);\mathbb L^2)}\\&= -\left( z(v-u^\tau),y(u^\tau,0)-y_d\right)_{L^2((0,T);\mathbb L^2)}-\left(u^\tau ,v-u^\tau\right)_{L^2(0,T)}\\&\leq \left(\norm{ z(v-u^\tau) }^2_{L^2((0,T); \mathbb L^2)}+\norm{v-u^\tau}^2_{\left(L^2(0,T)\right)^{N-1}}\right)^{1/2}\left(\norm{u^\tau}^2_{\left(L^2(0,T)\right)^{N-1}}+\norm{y(u^\tau,0)-y_d}^2_{L^2((0,T); \mathbb L^2)}\right)^{1/2}.
			\end{array}
		\end{equation*}
We can deduce  from \eqref{1inter} and \eqref{2inter} that
		\begin{equation}\label{5inter}
			\left(z(v-u^\tau), \frac{1}{\sqrt{\tau}}h^{ \tau}\right)_{L^2((0,T);\mathbb L^2)}\leq C \left(\norm{ z(v-u^\tau) }^2_{L^2((0,T); \mathbb L^2)}+\norm{v-u^\tau}^2_{\left(L^2(0,T)\right)^{N-1}}\right)^{1/2}.
		\end{equation}
		Set
		$$\mathcal{G}:=\left\{z(v):\; v\in\left(L^2(0,T)\right)^{N-1}\right\}.$$
		Then,  $\mathcal{G}\subset L^2((0,T); \mathbb L^2)$. We define on $\mathcal G\times \mathcal{ G}$ the inner product
		\begin{equation*}
			\langle z(v), z(w)\rangle_{\mathcal{G}}=(z(v),z(w))_{L^2((0,T);\mathbb L^2)}+(v,w)_{L^2(0,T)}.
		\end{equation*}
		Therefore, $\mathcal{G}$ endowed with the norm given by
		\begin{equation*}
			\norm{z(v)}^2_{\mathcal{G}}=\norm{ z(v) }^2_{L^2((0,T); \mathbb L^2)}+\norm{v}^2_{\left(L^2(0,T)\right)^{N-1}}\,\,\,\forall\, z(v)\in \mathcal{G},
		\end{equation*}
		is a Hilbert space.
		In view of \eqref{5inter}, we have that for every $v\in\left(L^2(0,T)\right)^{N-1}$,
		$$(z(v-u^\tau), \frac{1}{\sqrt{\tau}}h^{ \tau})_{L^2((0,T);\mathbb L^2)}\leq C \norm{z(v-u^\tau)}_{\mathcal{G}},
		$$ which implies that
		$$\norm{\frac{1}{\sqrt{\tau}}h^{ \tau}}_{\mathcal{G}^\star}\leq C.$$
		In particular, we have that
		\begin{equation}\label{6inter}
			\norm{\frac{1}{\sqrt{\tau}}h^{\tau}}_{L^2((0,T); \mathbb L^2)}\leq C.
		\end{equation}
		Moreover,  from \eqref{6inter}, \eqref{8Inter} and \eqref{paptau},  we can deduce that 
		$$\norm{p^\tau}_{L^2((0,T);\V)}\leq C.$$
Consequently, there exist $\zeta_1,\zeta_2\in L^2((0,T); \mathbb L^2)$ such that, as $\tau\to 0$,
		\begin{subequations}
			\begin{alignat}{6}
	\frac{1}{\sqrt{\tau}}\phi(u^\tau)&\rightharpoonup\,& \zeta_1&\,\text{weakly in} \,L^2((0,T); \mathbb L^2),\label{13inter}\\
	h^{\tau}&\rightharpoonup& h&\,\text{weakly in}\, L^2((0,T);\V)\label{9inter},\\\frac{1}{\sqrt{\tau}}h^{ \tau}&\rightharpoonup&\zeta_2&\,\text{weakly in}\, L^2((0,T); \mathbb L^2),\label{10inter}\\p^\tau&\rightharpoonup& p&\,\text{weakly in}\, L^2((0,T);\V)\label{11inter}.
			\end{alignat}
		\end{subequations}
Proceeding as for $y^{\tau}$ while using \eqref{13inter}-\eqref{11inter}, we get \eqref{paq} and \eqref{pap}. Thus,  we can deduce the identity  \eqref{7inter} while using \eqref{12inter}, \eqref{inter7} and \eqref{11inter}.  The proof is finished.
\end{proof}

\section{Concluding remarks}\label{conclusion}
We investigated an optimal control problem of a fractional parabolic partial differential equation involving a fractional Sturm-Liouville operator combined with the Caputo time-fractional derivative  in a space interval and in a general star graph; where the Sturm-Liouville operator is obtained as a composition of a left fractional Caputo derivative,  and a right fractional Riemann–Liouville derivative, and where the sources were missing. Using the no-regret control and the low regret control notions, we prove that the considered fractional optimal control admits a unique low regret control and we show that this control converges to the no-regret control that we characterize by a singular optimality system. Note that we can use the no- regret control theory for this kind of equation in the case were the boundary is missing, and we can as well use this technique to recover the missing information, this will be the object of a future investigation.

\bibliographystyle{abbrv}
\bibliography{references}

\end{document}